\documentclass[journal]{IEEEtran}
\newcommand{\beq}{\begin{equation}}
\newcommand{\eeq}{\end{equation}}
\newcommand{\bsq}{\begin{subequations}}
\newcommand{\esq}{\end{subequations}}
\newcommand{\bq}{\begin{eqnarray}}
\newcommand{\eq}{\end{eqnarray}}
\newcommand{\bqn}{\begin{eqnarray*}}
\newcommand{\eqn}{\end{eqnarray*}}
\usepackage{cite}
\usepackage{url}
\usepackage{empheq}
\usepackage{float}
\usepackage{dsfont}
\usepackage[hidelinks]{hyperref}

\usepackage{algorithm}
\usepackage{algorithmic}

\usepackage{amsmath,amssymb,amsfonts}
\allowdisplaybreaks[4]
\usepackage{graphicx}
\usepackage{textcomp}
\usepackage{amsmath}
\usepackage{enumerate}
\usepackage{epstopdf}
\usepackage{array}
\usepackage{booktabs}
\usepackage{subfigure}
\usepackage{multirow}
\usepackage{soul, color}
\usepackage[usenames,dvipsnames]{xcolor}

\newtheorem{theorem}{Theorem}

\newtheorem{corollary}{Corollary}
\newtheorem{proposition}{Proposition}

\soulregister\cite7 
\soulregister\citep7 
\soulregister\citet7 
\soulregister\ref7
\soulregister\pageref7

\makeatletter
\renewcommand{\maketag@@@}[1]{\hbox{\m@th\normalsize\normalfont#1}}%
\makeatother

\begin{document}

%
\title{Real-Time Bidding Strategy of Energy Storage in an Energy Market with Carbon Emission Allocation Based on Aumann-Shapley Prices}
%
%

\author{Rui Xie,~\IEEEmembership{Member,~IEEE},~Yue Chen,~\IEEEmembership{Member,~IEEE}
\thanks{This work was supported by the National Natural Science Foundation of China under Grant No. 52307144 and the Shun Hing Institute of Advanced Engineering, the Chinese University of Hong Kong, through Project RNE-p2-23 (Grant No. 8115071). (Corresponding to Y. Chen)}
\thanks{R. Xie and Y. Chen are with the Department of Mechanical and Automation Engineering, the Chinese University of Hong Kong, Hong Kong SAR. (email: ruixie@cuhk.edu.hk; yuechen@mae.cuhk.edu.hk)}
}

%
%

\markboth{Journal of \LaTeX\ Class Files,~Vol.~XX, No.~X, Feb.~2019}%
{Shell \MakeLowercase{\textit{et al.}}: Bare Demo of IEEEtran.cls for IEEE Journals}
%



\maketitle

\begin{abstract}
Energy storage (ES) can help decarbonize power systems by transferring green renewable energy across time. How to unlock the potential of ES in cutting carbon emissions by appropriate market incentives has become a crucial, albeit challenging, problem.
This paper fills the research gap by proposing a novel electricity market with carbon emission allocation and then investigating the real-time bidding strategy of ES in the proposed market. First, a carbon emission allocation mechanism based on Aumann-Shapley prices is developed and integrated into the electricity market clearing process to give combined electricity and emission prices. A parametric linear programming-based algorithm is proposed to calculate the carbon emission allocation more accurately and efficiently. Second, the real-time bidding strategy of ES in the proposed market is studied. To be specific, we derive the real-time optimal ES operation strategy as a function of the combined electricity and emission price using Lyapunov optimization. Based on this, the real-time bidding cost curve and bounds of ES in the proposed market can be deduced. Numerical experiments show the effectiveness and scalability of the proposed method. Its advantages over the existing methods are also demonstrated by comparisons.

\end{abstract}

\begin{IEEEkeywords}
carbon emission allocation, electricity-emission price, energy storage, real-time bidding, Lyapunov optimization
\end{IEEEkeywords}

%
\IEEEpeerreviewmaketitle

\section*{Nomenclature}
\addcontentsline{toc}{section}{Nomenclature}
\subsection{Abbreviations}
\begin{IEEEdescription}[\IEEEusemathlabelsep\IEEEsetlabelwidth{sssssss}]
\item[ES] Energy storage
\item[CEF] Carbon emission flow
\item[SoC] State of charge
\item[OPF] Optimal power flow
\item[LMP] Locational marginal price
\end{IEEEdescription}

       
\subsection{Parameters}
\begin{IEEEdescription}[\IEEEusemathlabelsep\IEEEsetlabelwidth{sssssss}]
\item[$\Psi_i$] Unit carbon emission of power plant $i$ 
\item[$D_{it}$] Load demand at bus $i$ in period $t$ 
\item[$L_i$] Loss sensitivity coefficient at bus $i$
\item[$L_0$] System loss linearization offset 
\item[$F_l$] Capacity of branch $l$ 
\item[$T_{li}$] Power transfer distribution factor between bus $i$ and branch $l$\item[$\underline{P}_{it}$/$\overline{P}_{it}$] Power output lower/upper bounds of power plant $i$ in period $t$ 
\item[$\tau$] Period length 
\item[$T$] Number of periods
\item[$\kappa$] Unit cost of carbon emission 
\item[$P_s^{max}$] Maximum power of ES $s$ 
\item[$\eta_s^c$/$\eta_s^d$] Charging/discharging efficiencies of ES $s$
\item[$\underline{E}_s$/$\overline{E}_s$] Lower/upper bounds of ES $s$'s stored energy
\item[$\underline{\gamma}_s$/$\overline{\gamma}_s$] Lower/upper bounds of the combined energy price of ES $s$ 
\item[$N_s$] Piecewise linearization parameter of the bidding cost curve of ES $s$
\end{IEEEdescription}

\subsection{Variables}
\begin{IEEEdescription}[\IEEEusemathlabelsep\IEEEsetlabelwidth{sssssss}]
\item[$p_{it}$] Net power output of power plant/ES $i$ in period $t$ 
\item[$f_{it}$] Energy cost of power plant/ES $i$ in period $t$ 
\item[$\lambda_{st}$/$\lambda_{it}$] LMPs at bus $s$/$i$ in period $t$ 
\item[$\lambda_t, \mu_{it}^\pm$] Dual variables of the OPF problem for electricity market clearing in period $t$
\item[$\mathcal{L}_t$] Lagrange dual function of the OPF problem for electricity market clearing in period $t$
\item[$\mathcal{E}$] Half of the total carbon emission cost of the power network 
\item[$\mathcal{E}_s$/$\mathcal{E}_i$] Carbon emission cost allocated to ES $s$/the load at bus $i$ 
\item[$\tilde{D}_i$] Combined power demand at bus $i$ considering the load and ES 
\item[$\psi_{st}$/$\psi_{it}$] Carbon emission prices at bus $s$/$i$ in period $t$ 
\item[$p_{st}^c$/$p_{st}^d$] Charging/discharging power of ES $s$ in period $t$ 
\item[$e_{st}$] ES $s$'s stored energy at the start of period $t$
\item[$\gamma_{st}$] Combined price of ES $s$ in period $t$ 
\item[$q_{st}$] Virtual queue of ES $s$ at the start of period $t$ 
\item[$E_s$] An adjustable parameter for establishing the virtual queue of ES $s$ 
\item[$l_{st}$] Lyapunov function of ES $s$ in period $t$ 
\item[$\Delta_{st}$] Lyapunov drift of ES $s$ in period $t$ 
\item[$V_s$] An adjustable parameter for Lyapunov optimization of ES $s$ 
\item[$\underline{P}_{st}$/$\overline{P}_{st}$] Power output lower/upper bounds of ES $s$ in period $t$ 
\end{IEEEdescription}

\section{Introduction}


\IEEEPARstart{M}{ore} than 100 countries have committed to achieving carbon neutrality in the 21st century to mitigate global climate change \cite{wei2022policy}. Reducing greenhouse gas emissions is one of the key steps. CO$_2$ is the main greenhouse gas whose emissions exceeded 37.9 Gt in 2021 \cite{united2022emissions}. Meanwhile, a large part of CO$_2$ emissions comes from fossil fuel-fired electricity generation \cite{kabeyi2022sustainable}. Therefore, low-carbon operation of power systems is a pivotal task toward the goal of carbon neutrality.

To promote low-carbon power system operation, an essential question is how to distribute carbon responsibilities among members in a power network. Carbon emissions are produced by fossil fuel power plants, but it is the consumers that create the electricity demand. In this regard, the consumers should be responsible for at least part of the carbon emissions \cite{kang2012carbon}. Since each consumer receives power from a mix of sources determined by Kirchhoff's law,
it is hard to quantify how many carbon emissions a consumer should be allocated. Moreover, energy storage (ES), a vital device for renewable integration, requires further consideration. 
While ES has a near-zero net energy consumption, it can help reduce power system carbon emissions by storing (releasing) electricity during periods with more (less) renewable energy. Hence, their carbon responsibilities should be allocated in a way that can maximize their potential for low-carbon power system operation.


One of the most commonly used methods to allocate carbon responsibilities is the carbon emission flow (CEF) method \cite{kang2015carbon}, which assumes that carbon emissions flow in the network along with the power flow. The CEF method has been applied in problems such as the operation scheduling \cite{wang2021optimal}, energy and carbon trading market \cite{lu2022peer}, and power system planning \cite{wei2021carbon}. The carbon emission allocation for ESs based on CEF was studied in \cite{wang2021optimal} and \cite{yang2023improved}. The ES was analogous to a container of liquid. CEF intensity and volume were used to describe the emissions related to the stored energy. However, the CEF method has some limitations:
1) the CEF result may change if virtual buses are added (an example is given in Section \ref{sec:CEF-limitation}); 2) the carbon intensity of ES under CEF only depends on its inflows during charging, but does not account for its outflows during discharging. This makes it hard for the CEF-based allocation to encourage ESs to shift more green energy into the periods with high carbon intensities.


Some other literature adopted cost-sharing mechanisms \cite{samet1982determination}. A comparison of those different mechanisms can be found in \cite{zhou2019cooperative}. It is revealed that the Shapley value performs the best but has high computational complexity, and the Aumann-Shapley pricing mechanism is a good alternative. In fact, the Aumann-Shapley pricing mechanism
is the unique cost-sharing mechanism characterized by scale invariance, consistency, additivity, and positivity axioms \cite{samet1982determination}. 
The Aumann-Shapley price-based carbon emission allocation mechanism was proposed in \cite{chen2018method} to promote power system emission reduction. This method was then used in low-carbon economic dispatch \cite{nan2022hierarchical} and carbon-trading-aware wind-battery planning \cite{nan2022bi}. However, the works above used numerical estimations of the partial derivatives and integrals when calculating the Aumann-Shapley value, which may be inaccurate or time-consuming. Moreover, ES was not considered.

The incorporation of ES is not trivial. To encourage ES to assist with low-carbon power system operation, we need to provide adequate carbon-oriented incentives (charges, prices) period by period, instead of allocating the carbon emissions in hindsight at the end of the day. Therefore, this paper aims to develop a real-time electricity market with carbon emission allocation and considers the participation of ES in the proposed market.
Literature related to real-time ES operation/bidding can be categorized into prediction-based and prediction-free ones \cite{wang2023online}. 
For the former category, model predictive control (MPC) was adopted in \cite{xie2021robust} to obtain the bidding strategy of wind-storage systems in a real-time market. It makes the current decisions based on the uncertainty predictions in future periods, which however are hard to obtain accurately in practice. 
Multi-stage stochastic/robust optimization is another prediction-based method. However, because they are computationally intractable, existing works solve them using simplification techniques such as affine policies, which scarifies optimality \cite{bodur2022two}.
In general, the prediction-based methods have their limitations in maintaining feasibility and pursuing optimality.

The Lyapunov optimization method developed by Neely \cite{neely2010stochastic} was originally used for stochastic network optimization in communication problems. This method can deal with stochastic programs with time average objectives and derive prediction-free online strategies, which make decisions based on the currently observed uncertainty realizations and do not depend on predictions. These advantages lead to its applications in online ES operation problems. 
A real-time coordinated ES and load operation strategy was proposed in \cite{li2016real}. An online and distributed strategy was introduced in \cite{zhong2019online} for shared ES. The online operation of a wind-ES system was studied in \cite{guo2021real}, whose optimal strategy was expressed using parametric linear programming techniques. The online energy management of microgrids consisting of electricity and heat generation and ESs was investigated in \cite{zhang2018online}. A real-time strategy for smart buildings equipped with ESs was proposed in \cite{ahmad2020real} considering the uncertainties from renewable generation, load demand, and energy prices. An online battery ES control algorithm was developed in \cite{shi2022lyapunov} to reduce peak loads and decrease electricity bills. However, the aforementioned works focused on the real-time operation of ES but did not address how ES bids in a real-time electricity market. The bidding problem is much more complicated than the operation problem since it needs to consider the impact of bids on the cleared price and quantity.
Moreover, none of these studies considered carbon emissions.

In this paper, an ES bidding cost curve is derived based on the online operation strategy by Lyapunov optimization, which can make the market clearing results the same as the online operation strategy under mild conditions. Compared with other methods for market participants’ optimal strategy, such as bilevel programming \cite{nasrolahpour2017bilevel}, the advantage of Lyapunov optimization is that it is prediction-free. Because of the development of renewable energy generation, the uncertainties in power systems have become more substantial. Moreover, the ES operator may lack information on the power network and other participants in the market, so it is difficult for the ES to predict future uncertainties and strategically participate in the market using bilevel programming. Instead, Lyapunov optimization does not rely on predictions and can provide an optimized strategy that only depends on the currently observed uncertainty realizations. Furthermore, derived from the Lyapunov optimization-based ES operation strategy, the ES bidding cost curve only depends on the ES’s current state-of-charge (SoC), so the ES operator no longer needs to know the current uncertainty realizations of renewable generation in the electricity market, which may not be accessible to the ES operator.

This paper fills the research gap by first proposing a novel electricity market with carbon emission allocation and then investigating the real-time bidding strategy of ES in the proposed market.
The contributions are as follows.

\emph{1) Electricity-emission Pricing.} In this paper, we propose an electricity market with carbon emission allocation. First, the electricity market is cleared by minimizing the total generation cost and total emission via lexicographic optimization. Then, the total emission is allocated among power plants, ESs, and loads based on Aumann-Shapley prices. A parametric linear programming-based algorithm is developed to calculate the emission prices. Compared to the existing Aumann-Shapley price calculation methods \cite{chen2018method,nan2022bi,nan2022hierarchical,zhou2019cooperative}, the proposed algorithm is faster and more accurate. The proposed carbon emission allocation method can avoid the limitations of the traditional CEF-based method \cite{kang2015carbon}.

\emph{2) Real-time ES Bidding Strategy.} To exploit the potential of ES in reducing carbon emissions by providing it with up-to-date carbon-oriented prices, the real-time bidding strategy of ES in the proposed market is studied. First, Lyapunov optimization is applied to establish the real-time optimal ES operation strategy as a function of the combined electricity and emission price. Compared to the existing work \cite{li2016real,zhong2019online,guo2021real,zhang2018online,ahmad2020real,shi2022lyapunov}, the proposed Lyapunov optimization method minimizes the exact drift-plus-penalty rather than its upper bound, which is shown to be more effective by the case studies. Then, the real-time ES bidding cost curve and bounds are derived. Numerical experiments show that the system total emission can be effectively reduced with ES in the proposed market.

The proposed carbon allocation mechanism can influence the market dynamics of wholesale electricity markets in several ways. First, it would impact the generation side. Carbon-intensive generators, like coal-fired plants, face higher operational costs, making them less competitive. In contrast, low-carbon alternatives such as renewables are encouraged. This can lead to increased investment in greener energy sources, thereby reducing carbon emissions. Second, it could affect the behavior of electricity consumers. When the emission is allocated to consumers according to their impacts on the system emission, the consumers may adjust their behavior, such as shifting demand to the period with more renewable generation. Third, the proposed carbon allocation mechanism can effectively quantify the contribution of ES to carbon emission reduction, and thus effectively incentivize the ES to help with low-carbon operation of the power system. Moreover, the proposed carbon allocation mechanism can also be combined with a cap-and-trade system \cite{carl2016tracking}, such as the EU Emissions Trading System (ETS) \cite{teixido2019impact} and the carbon trading policy in China \cite{zhang2020emission}. 

The rest of this paper is organized as follows. Section \ref{sec-2} introduces the proposed electricity market with carbon emission allocation. The real-time bidding strategy of ES is developed in Section \ref{sec-3}. The overall real-time market bidding and clearing procedures are also summarized. Case studies are presented in Section \ref{sec-5} with conclusions in Section \ref{sec-6}.

\section{Electricity Market with Emission Allocation}
\label{sec-2}

In this section, an electricity market with carbon emission allocation is proposed. 
We first clear the electricity market in Section \ref{sec2-a}. Then allocate the carbon emissions in Section \ref{sec2-b} and compare the proposed method with the CEF method in Section \ref{sec:CEF-limitation}.
For notation conciseness, the index $t$ of the time period is omitted in this section.

\subsection{Electricity Market Clearing}
\label{sec2-a}

In the proposed electricity market, the electricity prices are the locational marginal prices (LMPs) deduced from the Lagrangian function of an optimal power flow (OPF) problem. Particularly, we establish a lexicographic optimization-based OPF model, whose objective is
\begin{align}
\label{eq:lexicographic}
    \min_{p_i,\forall i}~ \left\{\sum_{i \in S_G \cup S_S} f_i (p_i), \sum_{i \in S_G} \sigma_i (p_i) \right\},
\end{align}
where $S_G$ and $S_S$ are the sets of power plants and ESs, respectively; $p_i$ is the power output of the power plant or ES $i$; $f_i(p_i)$ is the cost curve submitted to the market operator by the power plant or ES $i$, which is a piecewise linear and convex function of $p_i$ and has the unit \$/h; 
$\sigma_i(p_i)$ is the emission function of power plant $i$, which is also piecewise linear and convex. The two objective functions are optimized in lexicographic order, which means that the total generation cost $\sum_{i \in S_G \cup S_S} f_i (p_i)$ is minimized first, and then the total emission $\sum_{i \in S_G} \sigma_i(p_i)$ is minimized among all the feasible solutions with the minimum total power generation cost. 
In this way, the total emission is well-defined even if the OPF problem has multiple least-cost solutions. 

A lexicographic linear program can be equivalently transformed into a linear program with a weighted-sum objective function \cite{sherali1982equivalent}. Then, the objective \eqref{eq:lexicographic} can be replaced by
\begin{align}
    \min_{p_i, \forall i}~\sum_{i \in S_G \cup S_S} f_i (p_i) + \epsilon \sum_{i \in S_G} \sigma_i (p_i), \nonumber
\end{align}
for some small enough constant $\epsilon>0$. The optimal objective value is very close to the minimum generation cost.


The market clearing OPF problem is then formulated as
\begin{subequations}
\label{eq:energy-market}
\begin{align}
    \label{eq:energy-market-1}
    & \min_{p_i,\forall i} \!\sum_{i \in S_G \cup S_S} \!\!f_i (p_i) + \epsilon \sum_{i \in S_G} \sigma_i (p_i), \\
    \label{eq:energy-market-2} 
    & \mbox{s.t.} \sum_{i \in S_G \cup S_S} \!\!\!\!p_i \!-\! \sum_{i \in S_B}\!\! D_i =\!\!\!\!\! \sum_{i \in S_G \cup S_S} \!\!\!L_i p_i - \!\!\sum_{i \in S_B}\!\!\! L_i D_i + L_0: \bar \lambda, \\
    \label{eq:energy-market-3} 
    & \!\!-F_l\! \leq \!\!\!\!\sum_{i \in S_G \cup S_S}\!\!\!\! T_{li} p_i -\!\!\! \sum_{i \in S_B} T_{li} D_i \leq F_l: \mu_l^-, \mu_l^+ \geq 0, \forall l \in S_L, \\
    \label{eq:energy-market-4}
    & \underline{P}_i \leq p_i \leq \overline{P}_i, \forall i \in S_G \cup S_S,
\end{align}
\end{subequations}
where $S_B$ and $S_L$ are the index sets of buses and branches, respectively. $D_i$ is the load demand at bus $i$. $L_i$ is the loss sensitivity coefficient at bus $i$ and $L_0$ is the system loss offset, so the right side of \eqref{eq:energy-market-2} represents the system loss \cite{litvinov2004marginal} and \eqref{eq:energy-market-2} is the total power balance equation. The capacity of branch $l$ is denoted by $F_l$. $T_{li}$ is the power transfer distribution factor from bus $i$ to branch $l$.
Thus, \eqref{eq:energy-market-3} is the branch capacity constraint. The lower and upper bounds of the power output are stipulated in \eqref{eq:energy-market-4}. The bounds can be negative for ESs (How the ES bids the cost curve $f_i(p_i)$ and the bounds in real time will be discussed in Section \ref{sec-3}). 
Although the ramping limits of generators do not appear explicitly in the formulation of the proposed real-time market clearing problem, they can be considered by the generator owners during the bidding processes. The bid of a generator includes the cost curve $f_i (p_i)$ and the range of power output $[\underline{P}_i,\overline{P}_i]$. Since this is a real-time market \cite{litvinov2010design}, the power output range can be limited by the maximum and minimum outputs as well as the ramping requirements based on the previous output status.
The decision variables of the OPF problem \eqref{eq:energy-market} are $p_i, \forall i \in S_G \cup S_S$. $\overline{\lambda}$ and $\mu_l^\pm$ are dual variables.
Denote the Lagrangian function with constraints \eqref{eq:energy-market-2} and \eqref{eq:energy-market-3} by $\mathcal{L}$.
Then the LMP at bus $i$ is
\begin{align}
\label{eq:LMP}
    \lambda_i \triangleq \frac{\partial \mathcal{L}}{\partial D_i} = \bar \lambda (1 - L_i) + \sum_{l \in S_L} T_{li} (\mu_l^- - \mu_l^+), \forall i \in S_B.
\end{align}

To facilitate the computation of $\lambda_i$, considering that $f_i(p_i)$ and $\sigma_i(p_i)$ are piecewise linear and convex, they can be written as:
\begin{align}
    f_i (p_i) & = \max_{1 \leq n \leq N_i} \{ \alpha_{in}^F p_i + \beta_{in}^F \}, \forall i \in S_G \cup S_S, \nonumber \\
    \sigma_i (p_i ) & =\max_{1 \leq n \leq N_i} \{ \alpha_{in}^E p_i + \beta_{in}^E \}, \forall i \in S_G, \nonumber
\end{align}
where $\alpha_{in}^F, \beta_{in}^F, \alpha_{in}^E, \beta_{in}^E$ are parameters of each segment. 
Then, problem \eqref{eq:energy-market} is equivalent to the following linear program.
\begin{subequations}
\label{eq:market-LP}
\begin{align}
    & \min_{p_i,f_i,\sigma_i,\forall i}~ \sum_{i \in S_G \cup S_S} f_i + \epsilon \sum_{i \in S_G} \sigma_i, \\
    \label{eq:market-LP-2}
    & \mbox{s.t.}~ \eqref{eq:energy-market-2}-\eqref{eq:energy-market-4}, \\
    \label{eq:market-LP-3}
    & f_i \geq \alpha_{in}^F p_i + \beta_{in}^F,~ n = 1, \dots, N_i,~ \forall i \in S_G \cup S_S, \\
    \label{eq:market-LP-4}
    & \sigma_i \geq \alpha_{in}^E p_i + \beta_{in}^E,~ n = 1, \dots, N_i,~ \forall i \in S_G.
\end{align}
\end{subequations}

\emph{Remark:} Under the lexicographic order for generation cost and emission objective functions, the total cost is minimized first without emission objective, and then the total emission is minimized within the range where the total generation cost stays minimum. Therefore, if there is a unique optimal solution for the OPF without emission objective, this solution is also optimal for the lexicographic optimization.
We would like to clarify that instead of only considering the fuel cost, the first objective is to minimize the total generation cost, which is the summation of the bidding cost of power plants. Therefore, when power plants have to pay for their emissions, their bidding cost curve in the energy market will combine the fuel consumption and the emission costs.

There may be multiple OPF solutions under the single generation cost objective, which may lead to different total emissions. Therefore, we introduce the emission objective as the second objective, so that there is a unique total emission in the OPF solutions.
Nonetheless, the proposed emission allocation method is applicable to other types of objective functions. For instance, the weighted sum of electricity and emission cost was adopted in carbon-aware OPF literature \cite{lu2022peer}. Such an objective has the same form as the transformed objective function used in \eqref{eq:energy-market} and thus can be effectively handled by the proposed method.

The nodal loss coefficients depend on whether the node is injecting or consuming power, which can be addressed by a trial-and-error approach in the proposed energy market: First, initialize a power-flowing direction at each bus. Then solve the market clearing OPF problem. If the obtained direction is different from the assumed one, solve the OPF problem again under the updated direction. Repeat this process until the assumed and obtained directions are the same.

\subsection{Carbon Emission Allocation}
\label{sec2-b}
After the electricity market is cleared, we allocate the carbon emissions to the power plants, loads, and ESs. We introduce the proposed Aumann-Shapley price-based allocation mechanism, its properties, and the calculation algorithm.



\subsubsection{Allocation Mechanism}

Although carbon dioxide is only emitted in the power generation process, the demand side should take part of the responsibility because the load demands are the cause of power generation and the associated emissions. In the proposed carbon emission allocation mechanism, the power plants take responsibility for half of the emissions, and the other half is allocated to the ESs and loads. 
Specifically, power plant $i$ is responsible for $\sigma_i (p_i) \tau / 2$ emission in one period, where $\tau$ is the period length. 
Then, the cost curve and bounds of power plant $i$ are as follows.
\begin{align}
\label{eq:bidding-plant}
    f_i (p_i) = g_i (p_i) + \frac{1}{2} \kappa \cdot \sigma_i (p_i),~ \underline{P}_i \leq p_i \leq \overline{P}_i,~ \forall i \in S_G,
\end{align}
where $g_i(p_i)$ is the fuel cost function; $\kappa$ is the cost coefficient of emission and $\kappa \cdot \sigma_i (p_i) / 2$ is the emission cost function of power plant $i$. Equation \eqref{eq:bidding-plant} is the decomposition of the bidding cost curve, which is the same as the piecewise linear function mentioned in \eqref{eq:lexicographic} in Section \ref{sec2-a}.

In the following, we focus on how to allocate the other half of the total emission among ESs and load demands. We attribute carbon emissions to ESs because, despite consuming nearly zero energy across the whole time horizon, ESs have a positive or negative impact on carbon emissions in each period.
The key idea of the proposed Aumann-Shapley price-based allocation mechanism is that: we first treat the ES power outputs and load demands as given parameters, then derive how the total emission changes with them, and finally do the integral from 0 to their optimal strategies to allocate the emission costs. The detailed procedures are as follows:

First, given the ES power outputs $P_s,\forall s \in S_s$ and load demands $D_i,\forall i \in S_B$, we solve the following modified OPF problem:
\begin{align}
    & \min_{p_i,f_i,\sigma_i, \forall i}~ \sum_{i \in S_G \cup S_S} f_i + \epsilon \sum_{i \in S_G} \sigma_i, \nonumber \\
    \label{eq:market-fixed}
    & \mbox{s.t.}~ \eqref{eq:market-LP-2}-\eqref{eq:market-LP-4},~ p_s = P_s, \forall s \in S_S .
\end{align}
In the modified OPF problem \eqref{eq:market-fixed}, the ES power outputs $P=(P_s,\forall s \in S_S)$ and load demands $D=(D_i,\forall i \in S_B)$ are regarded as parameters. 
Denote the optimal solution of \eqref{eq:market-fixed} as $p_i^*,\forall i \in S_G$, which is a function of $P$ and $D$. With $p_i^*(P,D),\forall i \in S_G$, half of the total emission cost can be calculated by \eqref{eq:half-emission}, also a function of $P$ and $D$.
\begin{align}
\label{eq:half-emission}
    \mathcal{E}(P,D) \triangleq \frac{1}{2} \sum \nolimits_{i \in S_G} \kappa \cdot \sigma_i (p_i^*(P,D)) \cdot \tau.
\end{align}
Under given demands $D$, if $P_s$ is set as the optimal solution of $p_s,\forall s \in S_S$ in \eqref{eq:market-LP}, problem \eqref{eq:market-fixed} will produce the same optimal solution. To use Aumann-Shapley prices for emission allocation, we need the emission cost function $\mathcal{E}(P,D)$ under varying $P$ and $D$ (instead of merely the emission value under the optimal solution of \eqref{eq:market-LP}), which is defined in \eqref{eq:half-emission} using the optimal solution $p_i^*,\forall i \in S_G$ of problem \eqref{eq:market-fixed}.

The function $\mathcal{E}(P,D)$ can reflect how the system emission cost changes as the ES power outputs and load demands change. Then, the emission cost can be allocated as follows:
\begin{subequations}
\label{eq:emission-allocation}
\begin{align}
   \!\!\!\! & \mathcal{E}_s(P^*,D^*) \triangleq \!\!\int_0^{P_s^*} \!\!\frac{\partial \mathcal{E}}{\partial P_s} \left( \frac{y}{P_s^*}P^*, \frac{y}{P_s^*}D^* \right) \!dy, \forall s \in S_S, \\
    \!\!\!\! & \mathcal{E}_i(P^*,D^*) \triangleq \!\!\int_0^{D_i^*} \!\!\frac{\partial \mathcal{E}}{\partial D_i} \left(\frac{y}{D_i^*}P^*, \frac{y}{D_i^*}D^* \!\right) \!dy, \forall i \in S_B, 
\end{align}
\end{subequations}
where $y$ is the variable of the integration. It parametrizes the line segment from $(0,0)$ to $(P^*,D^* )$. 
$\mathcal{E}_s(P^*,D^*)$ and $\mathcal{E}_i(P^*,D^*)$ are the emission costs allocated to ES $s$ and load $i$, respectively, under ES power output $P^*$ and load demand $D^*$. In particular, $\mathcal{E}_s(P^*,D^*)$ is the integral of the partial derivative $\partial \mathcal{E}/\partial P_s$ along the segment from $(0,0)$ to $(P^*,D^*)$. $\partial \mathcal{E}/\partial P_s$ shows how the emission cost $\mathcal{E}$ changes as $P_s$ changes. The integral accumulates the influence of $P_s$ on $\mathcal{E}$. If ES $s$ discharges, it is likely to help decrease the total emission in that period and $\mathcal{E}_s(P^*,D^*)$ is negative. When $P_s^* = 0$, we have $\mathcal{E}_s = 0$. The function $\mathcal{E}(P,D)$ is determined by the modified OPF problem \eqref{eq:market-fixed}, whose result will not change when adding virtual buses. Hence, the proposed Aumann-Shapley price-based allocation mechanism can avoid the limitation of the CEF method (see Section \ref{sec:CEF-limitation}).



The emission price is the ratio of the allocated emission cost to the nonzero demand energy. The emission prices $\psi_s (P^*,D^*)$ of ES $s$ and $\psi_i (P^*,D^*)$ of load $i$ are
\begin{align}
\label{eq:emission-price}
    & \psi_s (P^*,D^*) \triangleq -\frac{1}{\tau} \int_0^1 \frac{\partial \mathcal{E}}{\partial P_s}(yP^*, yD^*)dy = \frac{\mathcal{E}_s(P^*, D^*)}{-P_s^* \tau}, \nonumber\\
    & \psi_i (P^*,D^*) \triangleq \frac{1}{\tau} \int_0^1 \frac{\partial \mathcal{E}}{\partial D_i}(yP^*, yD^*)dy = \frac{\mathcal{E}_i(P^*, D^*)}{D_i^* \tau}.
\end{align}
When $P_s^*$ or $D_i^*$ is $0$, $\psi_s$ and $\psi_i$ are still well-defined in \eqref{eq:emission-price}.

The emission costs are calculated when the ES charges or discharges, but the ES sometimes pays for it and sometimes gains income. In most cases, the emission price is positive. The charging of ES increases the requirement for power generation and causes system emission growth, so ES should pay for the emission under the proposed carbon emission allocation mechanism. In contrast, the discharging of ES generally decreases the system emission, and thus the ES receives revenue for its contribution to system emission reduction. In this way, the ES is encouraged to charge when the emission price is low (to pay less) and discharge when the emission price is high (to earn more), and then the potential of the ES in helping with system emission reduction is incentivized.

\subsubsection{Properties} 

Proposition \ref{prop:cost-sharing} below shows that the proposed method can ensure that half of the total emission is allocated to the ESs and loads.

\begin{proposition} 
\label{prop:cost-sharing}
Suppose problem \eqref{eq:market-fixed} is feasible for both $(0,0)$ and $(P^*,D^*)$, then the emission allocation in \eqref{eq:emission-allocation} is cost-sharing, i.e., \begin{align}
    \sum_{s \in S_S} \mathcal{E}_s(P^*, D^*) + \sum_{i \in S_B} \mathcal{E}_i(P^*, D^*) = \mathcal{E}(P^*, D^*) - \mathcal{E}(0, 0). \nonumber
\end{align}
\end{proposition}

Though the cost-sharing property has been proven for the Aumann-Shapley mechanism in \cite{samet1982determination}, it requires the total cost function to be continuously differentiable. In this paper, the function $\mathcal{E}(P,D)$ is not continuously differentiable, but we can still prove Proposition \ref{prop:cost-sharing} as in Appendix \ref{appendix-A}. Apart from the cost-sharing property, some other properties including scale invariance, monotonicity, additivity, and consistency are also listed in Appendix \ref{appendix-A}. 


Observe that in the modified OPF problem \eqref{eq:market-fixed}, the total emission depends on the net demand of loads and ESs at each bus. Therefore, we have the following corollary. 

\begin{corollary}
\label{cor:bus-price}
Every bus has a unique emission price, i.e., $\psi_s = \psi_i$ if $i = s \in S_S \subset S_B$, depending on the net demand $\tilde{D}_i$ at the bus, where
\begin{align}
    \tilde{D}_i \triangleq \left\{
    \begin{array}{ll}
        D_i, & \text{if}~ i \in S_B, i \notin S_S, \\
        D_i - P_s, & \text{if}~ i = s \in S_S \subset S_B.
    \end{array}
    \right.
    \label{eq:equivalent-load}
\end{align}
\end{corollary}




\emph{Remark:} In Proposition \ref{prop:cost-sharing}, the modified OPF problem \eqref{eq:market-fixed} is assumed to be feasible for both $(0,0)$ and $(P^*,D^*)$. Although the feasibility at $(P^*,D^*)$ is guaranteed by the fact that $P^*$ is part of the OPF solution under load $D^*$, the feasibility at $(0,0)$ may not hold due to the minimum power output constraints of thermal generators. To address this problem, we observe from the proof of Proposition 1 that if problem \eqref{eq:market-fixed} is feasible for $(P^0,D^0)$ and $(P^*,D^*)$, then
\begin{align}
    \mathcal{E}(P^*,D^*) - \mathcal{E}(P^0,D^0) = \sum\nolimits_{s \in S_S} I_s + \sum\nolimits_{i \in S_B} I_i, \nonumber
\end{align}
where
\begin{align}
    I_s \triangleq \int_{P_s^0}^{P_s^*} \frac{\partial \mathcal{E}}{\partial P_s} & \left( P^0 + \frac{y - P_s^0}{P_s^* - P_s^0} (P^* - P^0), \right. \nonumber \\
    & \left. D^0 + \frac{y - P_s^0}{P_s^* - P_s^0} (D^* - D^0) \right) dy,~ \forall s \in S_S, \nonumber \\
    I_i \triangleq \int_{D_i^0}^{D_i^*} \frac{\partial \mathcal{E}}{\partial D_i} & \left( P^0 + \frac{y - D_i^0}{D_i^* - D_i^0} (P^* - P^0), \right. \nonumber \\
    & \left. D^0 + \frac{y - D_i^0}{D_i^* - D_i^0} (D^* - D^0) \right) dy,~ \forall i \in S_B. \nonumber
\end{align}
If $(0,0)$ is an infeasible point, we can start from a feasible point $(P^0,D^0)$. Instead of using the definitions in \eqref{eq:emission-allocation}, the following emission cost allocation is adopted:
\begin{subequations}
\label{eq:emission-allocation-point}
\begin{align}
    & \mathcal{E}_s (P^*,D^*) = I_s + \mathcal{E}_s (P^0,D^0),~ \forall s \in S_S, \\
    & \mathcal{E}_i (P^*,D^*) = I_i + \mathcal{E}_i (P^0,D^0),~ \forall i \in S_B,
\end{align}
\end{subequations}
where
\begin{align}
    \sum\nolimits_{s \in S_S} \mathcal{E}_s (P^0,D^0) + \sum\nolimits_{i \in S_B} \mathcal{E}_i (P^0,D^0) = \mathcal{E} (P^0,D^0). \nonumber
\end{align}
Therefore, if there is a suitable method to allocate the emission cost at $(P^0,D^0)$, the emission cost increment from $(P^0,D^0)$ to $(P^*,D^* )$ can be allocated by Aumann-Shapley prices, and the allocation results in \eqref{eq:emission-allocation-point} still satisfy the cost-sharing property, i.e.,
\begin{align}
    \sum\nolimits_{s \in S_S} \mathcal{E}_s (P^*,D^*) + \sum\nolimits_{i \in S_B} \mathcal{E}_i (P^*,D^*) = \mathcal{E} (P^*,D^*). \nonumber
\end{align}
The remaining questions to answer are how to choose the feasible point $(P^0,D^0)$ and how to allocate the emission cost at $(P^0,D^0)$. We provide two ways in the following:

First, a day-ahead energy market clearing may provide a feasible point $(P^0,D^0)$ and the corresponding emission cost allocation. Since the day-ahead market can adjust unit commitment and consider the uncertainties by stochastic programming and robust optimization techniques, it can find a feasible point for real-time operation. It is worth noting that we do not need $P^0 \leq P^*$ or $D^0 \leq D^*$. In other words, the real-time adjustment can be either positive or negative whereas the cost-sharing property always holds. There are some existing works \cite{wang2022cross,ligao2023day,cheng2022carbon} on the carbon emission allocation in the day-ahead energy market. Generalizing the proposed method for the day-ahead energy market is also one of our future directions.

Second, $(P^0,D^0)$ can be chosen as the closest feasible point to $(0,0)$ on the line segment from $(0,0)$ to $(P^*,D^*)$, i.e., $(P^0,D^0)=(\zeta P^*,\zeta D^*)$, where $\zeta$ is the optimal value of the following linear program:
\begin{align}
    \min_{\zeta,p_i,D_i,\forall i}~ & \zeta \nonumber \\
    \mbox{s.t.}~ & \eqref{eq:energy-market-2}-\eqref{eq:energy-market-4},~ 0 \leq \zeta \leq 1, \nonumber \\
    & p_s = \zeta P_s^*,~ \forall s \in S_S,~ D_i = \zeta D_i^*,~ \forall i \in S_B. \nonumber
\end{align}
Then $(P^0,D^0)$ is feasible and $\mathcal{E}(P^0,D^0)$ can be calculated. We allocate the emission cost proportionally as follows:
\begin{align}
    \mathcal{E}_s (P^0,D^0) = \frac{-P_s^0 \cdot \mathcal{E} (P^0,D^0)}{-\sum_{s' \in S_S} P_{s'}^0 + \sum_{i \in S_B} D_i^0},~ \forall s \in S_S, \nonumber \\
    \mathcal{E}_i (P^0,D^0) = \frac{D_i^0 \cdot \mathcal{E}(P^0, D^0)}{- \sum_{s \in S_S} P_s^0 + \sum_{i' \in S_B} D_{i'}^0},~ \forall i \in S_B. \nonumber
\end{align}
Then the allocation results at $(P^*,D^*)$ satisfy the cost-sharing, scale invariance, monotonicity, additivity, and consistency properties, which can be directly verified using definitions. Therefore, even if the assumption of Proposition \ref{prop:cost-sharing} does not hold, we can slightly revise the allocation method and the calculation algorithm to maintain the desired properties.

\subsubsection{Calculation of Emission Prices}
We have seen that the proposed allocation mechanism possesses good properties. In the following, we develop an algorithm for calculating the emission prices efficiently and accurately.

For a fixed net demand vector $\tilde{D}$, the emission cost function \eqref{eq:half-emission} and the problem \eqref{eq:market-fixed} can be written in the following standard compact form.
\begin{align}
    \label{eq:compact}
   \mathcal{E}(\tilde{D}) = K^\top x~ \text{with}~ x ~\text{optimal in}~\! \left\{
\begin{aligned}
    \min_{x \geq 0}~ & C^\top x \\
    \mbox{s.t.}~ & A x = G \tilde{D} + H,
\end{aligned}
    \right.
\end{align}
where $x$ is a vector representing the decision variables in \eqref{eq:market-fixed}; $A$, $C$, $G$, $H$, and $K$ are constant matrices or vectors representing the coefficients.

Let $x^*$ be an optimal solution. According to linear programming theory \cite{dantzig2003linear}, we can divide $x^*$ into basic variables $x_B^*$ and nonbasic variables $x_N^*$ so that
\begin{align}
    x_B^* \geq 0,~ x_N^* = 0,~ x_B^* = A_B^{-1} (G \tilde{D} + H), \nonumber
\end{align}
with
\begin{align}
    x^* = \left(
    \begin{array}{cc}
        x_B^* \\
        x_N^*
    \end{array}
    \right),~
    A = \left( A_B, A_N \right),~
    K = \left(
    \begin{array}{cc}
        K_B \\
        K_N
    \end{array}
    \right), \nonumber
\end{align}
where $A_B$ is the optimal basis.

According to parametric linear programming theory \cite{dantzig2003linear}, when $\tilde{D}$ changes into $\tilde{D} + \Delta \tilde{D}$, $A_B$ remains the optimal basis if the basic variables are still nonnegative, i.e.,
\begin{align}
    A_B^{-1} (G (\tilde{D} + \Delta \tilde{D}) + H) \geq 0, \nonumber
\end{align}
and accordingly,
\begin{align}
    \mathcal{E}(\tilde{D} + \Delta \tilde{D}) & = K_B^\top A_B^{-1} (G (\tilde{D} + \Delta \tilde{D}) + H) \nonumber \\
    & = \mathcal{E}(\tilde{D}) + K_B^\top A_B^{-1} G \cdot \Delta \tilde{D}. \nonumber
\end{align}
Thus, we can calculate the partial derivative in \eqref{eq:emission-price} by
\begin{align}
\label{eq:partial-emission}
    \frac{\partial \mathcal{E}(\tilde{D})}{\partial \tilde{D}_i}  = \lim_{y \rightarrow 0} \frac{\mathcal{E}(\tilde{D} + y \omega_i) - \mathcal{E}(\tilde{D})}{y} = K_B^\top A_B^{-1} G \omega_i,
\end{align}
where $\omega_i$ is a constant vector with the same dimension as $\tilde{D}$ and it has $1$ at the $i$-th coordinate and $0$ elsewhere.

By the definition of emission price in \eqref{eq:emission-price}, we need to calculate the integral along a segment from $0$ to $\tilde{D}^*$. To do this, we start from $y_0 = 0$, determine the optimal basis $A_{B_0}$ at $(y_0+\delta)\tilde{D}^*$ ($\delta > 0$ is a small step length), and use \eqref{eq:partial-range} to find the interval $y \in [y_0,y_1]$ where $\partial \mathcal{E}/\partial \tilde{D}$ does not change.
\begin{align}
    A_{B_0}^{-1}(G \cdot y \tilde{D}^* + H) \geq 0,
    \label{eq:partial-range}
\end{align}
If $y_1 < 1$, we then move to $y_1+\delta$ and find a new interval $[y_1,y_2]$ by replacing $A_{B_0}$ with $A_{B_1}$, where $A_{B_1}$ is the optimal basis at $(y_1 + \delta)\tilde{D}^*$. Repeat until $y_M \geq 1$.


Let $y_M = 1$, then the emission price can be obtained by
\begin{align}
    \psi_i = \left[ \sum \nolimits_{m = 1}^M (y_m - y_{m-1}) K_{B_{m-1}}^\top A_{B_{m-1}}^{-1} G \omega_i\right]/\tau. \nonumber
\end{align}

The overall process is summarized in Algorithm \ref{alg:allocation}. Assume the segment from $0$ to $\tilde{D}^*$ meets $M'$ critical regions\footnote{A critical region is a region of parameters where the optimal basis does not change \cite{gal1972multiparametric}.}. By parametric linear programming theory \cite{gal1972multiparametric}, $M'$ is finite. 
With step parameter $\delta > 0$, the sample point number $m$ is no larger than $1/\delta$ and $M'$ when Algorithm \ref{alg:allocation} terminates. Therefore, Algorithm \ref{alg:allocation} terminates after finite and at most $O(\min \{ 1/\delta, M'\})$ steps. At each sample point, one linear OPF problem and three linear equation systems need to be solved. For a large-scale OPF problem and $M'>1/\delta$, the worst-case computational complexity is solving the OPF problem $1/\delta$ times. The result of Algorithm \ref{alg:allocation} is at least as accurate as the numerical calculation methods \cite{zhou2019cooperative,chen2018method,nan2022bi,nan2022hierarchical}. As $\delta$ decreases, the result becomes more accurate.
In addition, for a small enough $\delta > 0$, Algorithm \ref{alg:allocation} can encounter all the critical regions on the segment, so it gives the precise values of emission prices.

\begin{algorithm}
\normalsize
\caption{Emission Price Calculation}
\begin{algorithmic}[1]
\label{alg:allocation}
\renewcommand{\algorithmicrequire}{\textbf{Input:}}
\renewcommand{\algorithmicensure}{\textbf{Output:}}
\REQUIRE Parameters in \eqref{eq:market-fixed} and \eqref{eq:half-emission}; a step paramter $\delta > 0$. 
\ENSURE Emission prices $\psi_i, i \in S_B$.
\STATE Initiation: Calculate $A$, $C$, $G$, $H$, and $K$ in \eqref{eq:compact} and $\tilde{D}^*$ by \eqref{eq:equivalent-load}. Let $\psi_i \leftarrow 0, \forall i \in S_B$, $m \leftarrow 0$, $y_m \leftarrow 0$.
\STATE Let $m \leftarrow m + 1$. Solve the linear program in \eqref{eq:compact} with $\tilde{D} = (y_{m-1}+\delta) \tilde{D}^*$ to obtain the optimal basis $A_{B_{m-1}}$ and the corresponding $K_{B_{m-1}}$. Solve $A_{B_{m-1}}^{-1}(G \cdot y \tilde{D}^* + H) \geq 0$ and obtain an interval $[y', y'']$. Let $y_m \leftarrow \min\{y'',1\}$ and
$$    \psi_i \leftarrow \psi_i + \frac{1}{\tau}(y_m - y_{m-1}) K_{B_{m-1}}^\top A_{B_{m-1}}^{-1} G \omega_i,~ i \in S_B.$$
\STATE If $ y_m\geq 1$, terminate and output $\psi_i, i \in S_B$; otherwise, go to Step 2.
\end{algorithmic}
\end{algorithm} 

\subsection{Comparison with Carbon Emission Flow}
\label{sec:CEF-limitation}
Under the CEF method \cite{kang2015carbon}, carbon emissions are modeled as flows in the network along with the given power flows. 
For power flow $p_{ij}>0$ in directed line $(i,j)$, denote the corresponding CEF by $r_{ij}$. The CEF of the generation power $p_i$ at bus $i$ is denoted by $r_i^G$, which can be calculated according to the generation emission parameter. The CEF of the load demand $D_i$ is denoted by $r_i^D$. The CEF intensity $\rho$ is defined as the ratio of CEF to power. For example, $\rho_{ij}=r_{ij}/p_{ij}$. The CEF method has the following assumptions: For any $i \in S_B$ and $(i,j) \in S_L$,
\begin{subequations}
\begin{align}
    \label{eq:CEF-1}
    & r_i^G+ \sum_{(k,i) \in S_L} r_{ki} =r_i^D+ \sum_{(i,k) \in S_L} r_{ik}, \\
    \label{eq:CEF-2}
    & \rho_i^D = \rho_{ij} = \frac{r_i^D + \sum_{(i,k) \in S_L} r_{ik}}{D_i + \sum_{(i,k) \in S_L} p_{ik}} = \frac{r_i^G + \sum_{(k,i) \in S_L} r_{ki}}{p_i + \sum_{(k,i) \in S_L} p_{ki}},
\end{align}
\end{subequations}
where $S_B$ is the set of buses and $S_L$ is the set of directed lines with positive power flow; \eqref{eq:CEF-1} means the total emission flowing into a bus equals the total emission flowing out of it; \eqref{eq:CEF-2} indicates that the CEF flowing out of the same bus has the same carbon intensity. Based on the two assumptions, equations can be formulated under a given power flow distribution, whose solutions give the CEF values.

\begin{figure}[t]
\centering
\includegraphics[width=1.0\columnwidth]{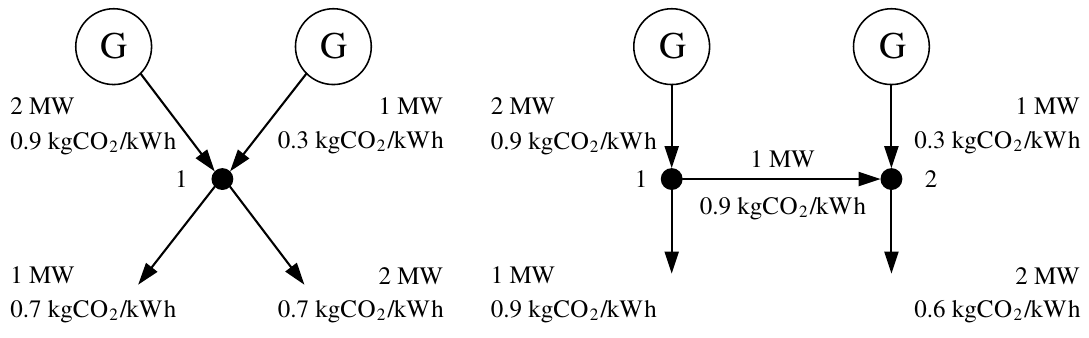}
\caption{An example to illustrate the shortcoming of the CEF method.}
\label{fig:CEFexample}
\end{figure}

However, the CEF method may have 
different allocation results if virtual buses are added. A simple example is depicted in Fig. \ref{fig:CEFexample}. There are two fossil fuel generators and two loads in the example system. 
On the left side of Fig. \ref{fig:CEFexample}, the carbon intensity flowing out of bus 1 is $(2 \times 0.9+1 \times 0.3)/(1+2)=0.7$ kgCO$_2$/kWh. On the right side, a virtual bus (bus 2) and a no-loss line are added to connect the two buses. Then the carbon intensity flowing out of bus 1 becomes $0.9$ kgCO$_2$/kWh, and that of bus 2 is $(1 \times 0.9+1 \times 0.3)/2=0.6$ kgCO$_2$/kWh.
The power flow outside the additional part remains the same as the original power flow. 

By removing all the virtual buses, the inconsistency of the CEF method can be avoided. However, the point of the example in Fig. \ref{fig:CEFexample} is that the allocation results of the CEF method may be controversial. For instance, if the line from bus 1 to bus 2 physically exists but is very short with a small impedance, it can be either modeled or not, without affecting the power flow calculation results much. In this circumstance, the load demand at bus 2 may prefer the model on the right side, where its carbon intensity is lower, but the load at bus 1 will be more favorable to the left model. There is still controversy due to the essential assumption of the CEF method that the power flowing out of a bus has the same carbon intensity.

The properties of the proposed Aumann-Shapley mechanism and the CEF method are summarized and compared in TABLE \ref{table:property}. The two methods both have cost-sharing and scale invariance properties. The monotonicity, additivity, and consistency properties depend on the emission cost function $\mathcal{E}(P,D)$. Since the CEF method allocates emissions according to the power flow rather than the function $\mathcal{E}(P,D)$, it does not have the three desired properties. When virtual buses are added to the network, the allocation results of the Aumann-Shapley mechanism will not change, while the CEF method can have different outcomes, as the example shows in Fig. 1. The proposed Aumann-Shapley mechanism can incentivize ES to help with system emission reduction, where the ES is rewarded for transferring green energy into high-emission periods. In contrast, the CEF method does not provide an effective incentive for ES, where the carbon intensity only depends on the inflows during charging but does not account for the discharging situation, so it cannot give a clear signal for the ES about when to discharge. The computational complexity of the CEF method involves solving a set of linear equations, which is very fast. We propose a calculation method to accelerate the computation of the Aumann-Shapley prices, which can reduce the computation burden down to solving some linear programs, comparable to the CEF method. To sum up, the Aumann-Shapley mechanism has more desired properties than the CEF method and outperforms the CEF method in the ability to encourage ES to help with system emission reduction. In addition, the computation speed is fast enough for a real-time market. Therefore, we recommend the proposed Aumann-Shapley mechanism for the carbon emission allocation in the energy market containing ESs.
\begin{table}[!ht]
\scriptsize
\renewcommand{\arraystretch}{1.3}
\caption{Comparison of Carbon Emission Allocation Methods}
\label{table:property}
\centering
\begin{tabular}{lcc}
\hline
& Aumann-Shapley & CEF \\
\hline
Cost-sharing & \checkmark & \checkmark \\
Scale invariance & \checkmark & \checkmark \\
Monotonicity & \checkmark & - \\
Additivity & \checkmark & - \\
Consistency & \checkmark & - \\
Invariance after adding virtual buses & \checkmark & $\times$ \\
Incentivizing ES & \checkmark & $\times$ \\
Computation speed & Fast & Very fast \\
\hline
\end{tabular}
\end{table}

\section{Real-time Energy Storage Bidding Strategy}
\label{sec-3}
In Section \ref{sec-2}, we develop an electricity market with carbon emission allocation, then a remaining question is how the ES $s$ determines its real-time bidding cost curve $f_s(p_s)$ and bounds $\underline{P}_s, \overline{P}_s$. The design of a bidding cost curve follows the rule that the resulting optimal dispatch strategy of ES $s$ by \eqref{eq:energy-market} should be the same as the optimal operation strategy of ES $s$ under the corresponding combined electricity and emission price $\lambda_{st}+\psi_{st}$ \cite{litvinov2010design}. 
In this section, we first assume the energy and emission prices are given and analyze the optimal ES operation strategy by developing an offline model and its online counterpart based on Lyapunov optimization. After the relationship between the real-time optimal ES operation strategy and the price is obtained, the bidding cost curve and bounds are derived.


\subsection{Offline Energy Storage Operation Model}

Suppose the energy prices $\lambda_{st},\forall t=1,\dots,T$ and the emission prices $\psi_{st},\forall t=1,\dots,T$ are given. In the offline optimal operation model of ES $s$, the charging power $p_{st}^c$ and discharging power $p_{st}^d$ in each period are optimized to maximize the total revenue as follows.
\begin{subequations}
\label{eq:offline}
\begin{align}
    \label{eq:offline-1} 
    & \textbf{P1}: ~\max_{p_{st}^c, p_{st}^d, e_{st},\forall t} ~\sum \nolimits_{t = 1}^T (\lambda_{st} + \psi_{st}) (p_{st}^d - p_{st}^c) \tau, \\
    \label{eq:offline-2} 
   &  \mbox{s.t.} ~ 0 \leq p_{st}^c \leq P_s^{max}, 0 \leq p_{st}^d \leq P_s^{max}, p_{st}^c p_{st}^d = 0, \forall t, \\
    \label{eq:offline-3} 
    & e_{s(t+1)} = e_{st} + p_{st}^c \tau \eta_s^c - p_{st}^d \tau / \eta_s^d, \forall t, \\
    \label{eq:offline-4}
    & \underline{E}_s \leq e_{st} \leq \overline{E}_s, \forall t, 
\end{align}
\end{subequations}
where the maximized objective function in \eqref{eq:offline-1} is the total revenue considering both energy and emission. Constraint \eqref{eq:offline-2} stipulates bounds for charging and discharging power and prohibits simultaneous charging and discharging. $e_{st}$ denotes the stored energy at the beginning of period $t$. $\eta_s^c$ and $\eta_s^d$ are the charging and discharging efficiencies, respectively. Constraint \eqref{eq:offline-3} reflects the state-of-charge (SoC) dynamics and \eqref{eq:offline-4} sets lower and upper bounds for the stored energy.

The offline model \eqref{eq:offline} can help the ES get the maximum revenue over the whole time horizon. However, the ES cannot derive its real-time bidding strategy directly based on the offline model because: 1) The optimal solution of \eqref{eq:offline} depends on the combined electricity and emission prices for all periods, while what we need in a real-time market is a bidding cost curve depending solely on the current period. 2) The model \eqref{eq:offline} fails to reflect the influence of the ES charging and discharging strategies on future electricity and emission prices. To overcome the above limitations, we propose a real-time optimal ES operation model to derive the real-time ES bidding strategy in the following.



\subsection{Real-Time Optimal Energy Storage Operation Strategy}

The offline problem \eqref{eq:offline} maximizes the total revenue in $T$ periods. Since actually we want to maximize the long-term time average revenue, let $T \rightarrow \infty$ and transform \eqref{eq:offline} into the following problem.
\begin{align}
    & \textbf{P1}':~v_0^* \triangleq - \min_{p_{st}^c, p_{st}^d, e_{st}, \forall t}~  \lim_{T \rightarrow \infty} \frac{1}{T} \sum_{t = 1}^T \mathbb{E} [-\gamma_{st} (p_{st}^d - p_{st}^c) \tau], \nonumber \\
    \label{eq:v_0}
    & \mbox{s.t.}~ \eqref{eq:offline-2}-\eqref{eq:offline-4}. 
\end{align}
In \eqref{eq:v_0}, the combined energy and emission price $\gamma_{st}$ is defined by $\gamma_{st} \triangleq \lambda_{st} + \psi_{st}$, which is the source of uncertainty in the optimization; $\mathbb{E}(\cdot)$ means taking an expectation; maximizing the average expected revenue is equivalently transformed into minimizing its negative value; the optimal time average expected revenue is denoted by $v_0^*$.

The problem \eqref{eq:v_0} is still an offline model. Then, we use Lyapunov optimization to turn \eqref{eq:v_0} into its online counterpart. First, we define a virtual queue $q_{st}, t = 1, 2, \dots$ as follows.
\begin{align}
    q_{st} \triangleq e_{st} - E_s, \forall t, \nonumber
\end{align}
where $E_s$ is a parameter to be determined later. By \eqref{eq:offline-4}, $\underline{E}_s - E_s \leq q_{st} \leq \overline{E}_s - E_s, \forall t$, so $\lim_{T \rightarrow \infty} \mathbb{E}[|q_{sT}|]/T = 0$, which means the virtual queue $\{q_{st},\forall t\}$ is mean rate stable \cite{neely2010stochastic}. Then we relax constraint \eqref{eq:offline-4} to the mean rate stability of the virtual queue $\{q_{st},\forall t\}$, restate constraint \eqref{eq:offline-3} using $q_{st}$, and obtain a relaxed problem as follows.
\begin{align}
    &  v_1^* \triangleq - \min_{p_{st}^c, p_{st}^d, q_{st},\forall t}~ \lim_{T \rightarrow \infty} \frac{1}{T} \sum_{t = 1}^T \mathbb{E} [-\gamma_{st} (p_{st}^d - p_{st}^c) \tau], \nonumber \\
     & \mbox{s.t.}~ \eqref{eq:offline-2},~ \lim_{T \rightarrow \infty} \frac{1}{T} \mathbb{E} [|q_{sT}|] = 0, \nonumber \\
    \label{eq:v_1}
   & q_{s(t+1)} = q_{st} + p_{st}^c \tau \eta_s^c - p_{st}^d \tau / \eta_s^d, \forall t .
\end{align}

Define a Lyapunov function $l_{st}$ for the virtual queue $q_{st}$ by $l_{st} \triangleq (q_{st})^2/2$. Define Lyapunov drift $\Delta_{st}$ as the change of the Lyapunov function:
\begin{align}
    \Delta_{st} & \triangleq l_{s(t+1)}-l_{st} = (q_{s(t+1)})^2/2 - (q_{st})^2/2 \nonumber \\
    & = (p_{st}^c \tau \eta_s^c - p_{st}^d \tau / \eta_s^d)^2/2 + (p_{st}^c \tau \eta_s^c - p_{st}^d \tau / \eta_s^d)q_{st}, \forall t. \nonumber
\end{align}

Then, we can derive an online algorithm by minimizing the weighted sum of the Lyapunov drift $\Delta_{st}$ and the objective for that specific period $t$ as follows.
\begin{align}
   & \textbf{P2}:~\min_{p_{st}^c, p_{st}^d}~ \Delta_{st} + V_s(-\gamma_{st} (p_{st}^d - p_{st}^c) \tau), \nonumber \\
   & \mbox{s.t.}~ \Delta_{st} = (p_{st}^c \tau \eta_s^c - p_{st}^d \tau / \eta_s^d)^2/2 + (p_{st}^c \tau \eta_s^c - p_{st}^d \tau / \eta_s^d)q_{st}, \nonumber \\
   \label{eq:drift-plus-penalty}
   & 0 \leq p_{st}^c \leq P_s^{max}, 0 \leq p_{st}^d \leq P_s^{max}, p_{st}^c p_{st}^d = 0,
\end{align}
where $V_s > 0$ is a penalty coefficient to be determined later and the objective function is called a drift-plus-penalty term. By minimizing the drift-plus-penalty term, we can balance the virtual queue stability and the time average expectation of revenue. It is worth noting that distinct from the previous works \cite{li2016real,zhong2019online,guo2021real,zhang2018online,ahmad2020real,shi2022lyapunov} that used an upper bound of $\Delta_{st}$ and minimized the linear function $C_0+(p_{st}^c \tau \eta_s^c - p_{st}^d \tau / \eta_s^d)q_{st}+V_s(-\gamma_{st} (p_{st}^d - p_{st}^c) \tau)$ with a constant $C_0$, we minimize the exact drift-plus-penalty term, which is a quadratic function. This can improve the accuracy of the online algorithm.

In each period $t$, based on the up-to-date $q_{st}$ and $\gamma_{st}$, problem \eqref{eq:drift-plus-penalty} can be solved to obtain $p_{st}^c$, $p_{st}^d$, and $q_{s(t+1)}$. With the condition $p_{st}^c p_{st}^d = 0$ and the expression of $\Delta_{st}$, problem \eqref{eq:drift-plus-penalty} is equivalently transformed into
\begin{align}
    \min \left\{ \min_{0 \leq p_{st}^c \leq P_s^{max}} \left\{ (p_{st}^c \tau\eta_s^c )^2 / 2 + p_{st}^c \tau q_{st} \eta_s^c + V_s \gamma_{st} p_{st}^c \tau \right\}, \right. \nonumber \\
    \left. \min_{0 \leq p_{st}^d \leq P_s^{max}} \left\{ (p_{st}^d \tau / \eta_s^d)^2 / 2 - p_{st}^d \tau q_{st} / \eta_s^d - V_s \gamma_{st} p_{st}^d \tau \right\}
    \right\}, \nonumber
\end{align}
which comes down to finding the minimum of constrained 1-dimensional quadratic functions. Define the net output of ES $s$ by $p_{st} \triangleq p_{st}^d - p_{st}^c$, which is consistent with our notations in the proposed electricity market. Then the optimal solution $p_{st}$ of problem \eqref{eq:drift-plus-penalty} is a piecewise linear function of the combined electricity and emission price $\gamma_{st}$:
\begin{align}
\small
\label{eq:strategy-gamma}
    \left\{
    \begin{array}{ll}
        - P_s^{max}, & \text{if}~ q_{st} \leq - \frac{V_s \gamma_{st}}{\eta_s^c} - P_s^{max} \tau\eta_s^c , \\
        \frac{q_{st} \eta_s^c + V_s \gamma_{st}}{\tau(\eta_s^c)^2 }, & \text{if}~ - \frac{V_s \gamma_{st}}{\eta_s^c} - P_s^{max} \tau \eta_s^c \leq q_{st} \leq - \frac{V_s \gamma_{st}}{\eta_s^c}, \\
        0, & \text{if}~ - \frac{V_s \gamma_{st}}{\eta_s^c} \leq q_{st} \leq - V_s \gamma_{st} \eta_s^d, \\
        \frac{q_{st} / \eta_s^d + V_s \gamma_{st}}{\tau / (\eta_s^d)^2}, & \text{if}~ - V_s \gamma_{st} \eta_s^d \leq q_{st} \leq - V_s \gamma_{st} \eta_s^d + \frac{P_s^{max} \tau}{\eta_s^d}, \\
        P_s^{max}, & \text{if}~ - V_s \gamma_{st} \eta_s^d + \frac{P_s^{max} \tau}{\eta_s^d} \leq q_{st}.
    \end{array}
    \right.
\end{align}
With the online strategy \eqref{eq:strategy-gamma}, we can derive the real-time bidding strategy in Section \ref{sec:bidding}. Before that, we first determine the values of parameters $E_s$ and $V_s$ to guarantee the SoC range constraint \eqref{eq:offline-4}, which is relaxed in \textbf{P2}. $E_s$ and $V_s$ can be set according to the theorems below.


\begin{theorem}
\label{thm:feasible}
    Assume $\gamma_{st} \in [\underline{\gamma}_s, \overline{\gamma}_s]$ with $0 \leq \underline{\gamma}_s < \overline{\gamma}_s \eta_s^c \eta_s^d$. If $E_s$ and $V_s$ satisfy
    \begin{subequations}
    \label{eq:V-E-range}
    \begin{align}
        \label{eq:V-E-range-1}
        & 0 < V_s \leq \frac{\eta_s^c (\overline{E}_s - \underline{E}_s)}{\overline{\gamma}_s \eta_s^c \eta_s^d - \underline{\gamma}_s},  \\
        \label{eq:V-E-range-2}
        & \underline{E}_s + V_s \overline{\gamma}_s \eta_s^d \leq E_s \leq \overline{E}_s + V_s \underline{\gamma}_s / \eta_s^c, 
    \end{align}
    \end{subequations}
    the SoC range constraint $\underline{E}_s \leq e_{st} \leq \overline{E}_s, \forall t$ holds automatically under the operation strategy in \eqref{eq:strategy-gamma}. 
\end{theorem}

The proof of Theorem \ref{thm:feasible} is in Appendix \ref{appendix-B}. In the assumption of Theorem \ref{thm:feasible}, $[\underline{\gamma}_s, \overline{\gamma}_s]$ is the range of the combined price $\gamma_{st}, \forall t$ that can be estimated using historical data. 
The assumption $\underline{\gamma}_s < \overline{\gamma}_s \eta_s^c \eta_s^d$ is reasonable as it is the condition for ES $s$ to possibly make positive profits considering the loss in the charging and discharging process. 

Another issue we care about is the gap between the online result by \textbf{P2} and that of $\textbf{P1}'$, as discussed below.

\begin{theorem}
\label{thm:performance}
    Let the parameters $V_s$ and $E_s$ be in the range of \eqref{eq:V-E-range}. Assume $\gamma_{st}, \forall t$ are independent and identically distributed. Denote the time average revenue expectation of the strategy \eqref{eq:strategy-gamma} (which is also the optimal solution of \textbf{P2}) by $v^*$, then 
    \begin{align}
        v_0^*-(P_s^{max} \tau)^2 / (2V_s(\eta_s^d)^2) \leq v^* \leq v_0^*. \nonumber
    \end{align}
\end{theorem}

The proof of Theorem \ref{thm:performance} is in Appendix \ref{appendix-C}. On the one hand, a larger $V_s$ leads to a tighter performance bound according to Theorem \ref{thm:performance}. On the other hand, Theorem \ref{thm:feasible} limits the choice of parameters $E_s$ and $V_s$. To achieve the best online performance, we maximize $V_s$ concerning the constraints in \eqref{eq:V-E-range}. Then, $E_s$ and $V_s$ are chosen as
\begin{align}
\label{eq:V-E-value}
    E_s = \frac{\overline{\gamma}_s \eta_s^c \eta_s^d \overline{E}_s - \underline{\gamma}_s \underline{E}_s}{\overline{\gamma}_s \eta_s^c \eta_s^d - \underline{\gamma}_s},~ V_s = \frac{\eta_s^c (\overline{E}_s - \underline{E}_s)}{\overline{\gamma}_s \eta_s^c \eta_s^d - \underline{\gamma}_s}.
\end{align}

\emph{Remark:} We adopt the assumption $\underline{\gamma}_s \geq 0$ for theoretical convenience considering the fact that negative price rarely arises \cite{seel2021plentiful}. The performance guarantee and parameter choices derived in Theorem \ref{thm:feasible} and Theorem \ref{thm:performance} can serve as a reference for the general cases with possibly negative prices. It is worth noting that the feasibility of the bidding strategy to be proposed in \eqref{eq:bidding-ES} does not depend on the assumption of positive prices, which will be elaborated in Proposition \ref{prop:bidding}.

\subsection{Real-time Bidding Strategy of Energy Storage}
\label{sec:bidding}

In the following, we derive the real-time bidding strategy of ES based on \eqref{eq:strategy-gamma}. 
First, with the proposed real-time optimal operation strategy \eqref{eq:strategy-gamma}, the net output of ES $s$ is a function $p_{st}(\gamma_{st})$ of the combined price $\gamma_{st}$. The function $p_{st}(\gamma_{st})$ is nondecreasing and its minimum $\underline{P}_{st}$ and maximum $\overline{P}_{st}$ are

\begin{footnotesize}
\begin{align}
    & \underline{P}_{st} = \left\{
    \begin{array}{ll}
        - P_s^{max}, & \text{if}~ q_{st} \leq - \frac{V_s \underline{\gamma}_s}{\eta_s^c} - P_s^{max} \tau \eta_s^c, \\
        \frac{q_{st} \eta_s^c + V_s \underline{\gamma}_s}{\tau(\eta_s^c)^2 }, & \text{if}~ - \frac{V_s \underline{\gamma}_s}{\eta_s^c} - P_s^{max} \tau \eta_s^c \leq q_{st} \leq - \frac{V_s \underline{\gamma}_s}{\eta_s^c}, \\
        0, & \text{if}~ q_{st} \geq - \frac{V_s \underline{\gamma}_s}{\eta_s^c}.
    \end{array}
    \right. \nonumber \\
    \label{eq:bounds-ES}
    & \overline{P}_{st} = \left\{
    \begin{array}{ll}
        0, \!\!\!\!& \text{if}~ q_{st} \leq - V_s \overline{\gamma}_s \eta_s^d, \\
        \frac{q_{st} / \eta_s^d + V_s \overline{\gamma}_s}{\tau / (\eta_s^d)^2}, \!\!\!\!& \text{if}~ - V_s \overline{\gamma}_s \eta_s^d \leq q_{st} \leq - V_s \overline{\gamma}_s \eta_s^d + \frac{P_s^{max} \tau}{\eta_s^d}, \\
        P_s^{max}, \!\!\!\!& \text{if}~ - V_s \overline{\gamma}_s \eta_s^d + \frac{P_s^{max} \tau}{\eta_s^d} \leq q_{st}.
    \end{array}
    \right. 
\end{align}
\end{footnotesize}

The idea of creating ES bidding curves is that the resulting dispatch strategy $p_{st}$ of ES $s$ in the market clearing should be the same as the optimal operation strategy $p_{st}(\gamma_{st})$ of ES under the combined electricity and emission price $\gamma_{st} = \lambda_{st} + \psi_{st}$ \cite{litvinov2010design}. The market clearing OPF problem \eqref{eq:energy-market} minimizes the total cost, aiming at seeking the social optimum.
According to strong duality, the optimal solution to  \eqref{eq:energy-market} can be equivalently obtained by $\max_{\lambda_{t}} \min_{p \in [\underline{P},\overline{P}]} \mathcal{L}_t$. Given the LMP $\lambda_t^*$, the inner minimization problem can be separated into the problems for each agent. In particular, for ES $s$, it solves
\begin{align} 
    & \min_{p_{st} \in [\underline{P}_{st}, \overline{P}_{st}]} f_{st}(p_{st}) - \lambda_{st}^* p_{st} \nonumber \\
    \label{eq:ESindividual}
    =~ &  \min_{p_{st} \in [\underline{P}_{st}, \overline{P}_{st}]} f_{st}(p_{st}) - \gamma_{st}^* p_{st} + \psi_{st}^* p_{st},
\end{align}
where $\psi_{st}^*$ and $\gamma_{st}^*$ are the emission price and the combined electricity and emission price in the market clearing and carbon emission allocation results. To implement the aforementioned idea, the desired bidding curve $f_{st}(p_{st})$ should be convex in $p_{st}$ and make the optimal operation strategy $p_{st}(\gamma_{st})$ become the optimal solution of \eqref{eq:ESindividual}. Therefore, we choose the following bidding cost curve:
\begin{align}
    f_{st}(p_{st}) \triangleq \int_0^{p_{st}} \gamma_{st} dp_{st} - \psi_{st}^* p_{st},~ \underline{P}_{st} \leq p_{st} \leq \overline{P}_{st}. \nonumber
\end{align}
Using the formula of the optimal operation strategy $p_{st}(\gamma_{st})$ in \eqref{eq:strategy-gamma}, we can get $f_{st}(p_{st}) = $
\begin{align}
\left\{
\begin{array}{ll}
\frac{p_{st} \eta_s^c (p_{st} \tau \eta_s^c - 2 q_{st})}{2V_s} - \psi_{st}^* p_{st}, & \text{if}~ \underline{P}_{st} \leq p_{st} \leq 0, \\
\frac{p_{st} (p_{st} \tau - 2 q_{st} \eta_s^d)}{2 V_s (\eta_s^d)^2} - \psi_{st}^* p_{st}, & \text{if}~ 0 \leq p_{st} \leq \overline{P}_{st}.
\end{array}
\right.
\label{eq:cost-function}
\end{align}
The effectiveness of the function $f_{st}(p_{st})$ is demonstrated in Proposition \ref{prop:bidding} with proofs in Appendix \ref{appendix-D}.

\begin{proposition}
\label{prop:bidding} 
The function $f_{st}(p_{st})$ is convex. Moreover, when $f_{st}(p_{st})$ with $\underline{P}_{st} \leq p_{st} \leq \overline{P}_{st}$ is used as the bidding cost curve, the operational constraints of ES must hold (even if $\gamma_{st} \notin [\underline{\gamma}_s, \overline{\gamma}_s]$) and the dispatch strategy $p_{st}$ determined by the market clearing problem \eqref{eq:energy-market} coincides with the operation strategy in \eqref{eq:strategy-gamma}.
\end{proposition}

It is worth noting that no assumption is needed in Proposition \ref{prop:bidding}. When $\kappa = 0$, the carbon emissions cause no costs, and the combined price degenerates to the electricity price, which is the same as the situation in a traditional LMP-based electricity market \cite{wu2013impact}. However, the definition of $f_{st}(p_{st})$ in \eqref{eq:cost-function} is not practical, because the ES operator cannot know the emission price $\psi_{st}^*$ in advance. Therefore, we approximate $\psi_{st}^*$ by the known value $\psi_{s(t-1)}$ from the previous period, i.e., $f_{st}(p_{st}) \approx$
\begin{align}
\left\{
\begin{array}{ll}
\frac{p_{st} \eta_s^c (p_{st} \tau \eta_s^c - 2 q_{st})}{2V_s} - \psi_{s(t-1)} p_{st}, & \text{if}~ \underline{P}_{st} \leq p_{st} \leq 0, \\
\frac{p_{st} (p_{st} \tau - 2 q_{st} \eta_s^d)}{2 V_s (\eta_s^d)^2} - \psi_{s(t-1)} p_{st}, & \text{if}~ 0 \leq p_{st} \leq \overline{P}_{st},
\end{array}
\right. \nonumber
\end{align}
which is convex.

Since in the proposed electricity market, the cost curves for power plants are piecewise linear, to be consistent, we approximate $f_{st}(p_{st})$ by a piecewise linear function: Choose $N_s$ points from the interval $[\underline{P}_{st}, \overline{P}_{st}]$ and calculate the function values as follows.
\begin{align}
    & \underline{P}_{st} = P_{st1} < \dots < P_{stn} < \dots < P_{stN_s} = \overline{P}_{st}, \nonumber \\
    & F_{stn} \triangleq g_{st}(P_{stn})- \psi_{s(t-1)} P_{stn},~ n = 1, 2, \dots, N_s. \nonumber
\end{align}
Then for any $p_{st} \in [\underline{P}_{st}, \overline{P}_{st}]$,
\begin{align}
& f_{st}(p_{st}) \approx \nonumber \\
\label{eq:bidding-ES}
& \max_{1 \leq n \leq N_s-1} \left\{ F_{stn} + \frac{F_{st(n+1)} - F_{stn}}{P_{st(n+1)} - P_{stn}} (p_{st} - P_{stn}) \right\}, 
\end{align}
where the right side is the final cost curve submitted by ES $s$.

\emph{Remark:} For the dispatch problem in a power system with ESs, the operation of the ES is optimized to minimize the total cost of the system. The system operator has the information of all the components of the system, and possibly has valuable predictions for future uncertainties such as renewable generation and load demand, which helps to approach social optimum via stochastic programming. However, in the considered scenarios, the ES participates in the energy market for its profits. Meanwhile, the ES does not have the information of other participants in the market, nor does the ES have access to the predictions of uncertainties. Therefore, the ES submits the bidding cost curve according to the current SoC value.

Traditional real-time ES operation methods \cite{li2016real,zhong2019online,guo2021real,zhang2018online,ahmad2020real,shi2022lyapunov} did not address how ES bids in the market. In contrast, the impacts of the ES bidding on the LMP are considered in this paper. The ES operator creates the ES bidding curves to make profits. However, the ES operator is not a price-maker in the proposed method. Otherwise, the market clearing problem should be considered as a lower-level problem in the bilevel bidding curve optimization problem \cite{nasrolahpour2017bilevel}, which is much more complicated when the ES operator does not have enough information on the power network and other market participants, such as future renewable generation. Therefore, we leave the strategic behavior of price-making ES operators for future work.

\subsection{Overall Procedure}
\label{sec-4}

\begin{figure}[t]
\centering
\includegraphics[width=1.0\columnwidth]{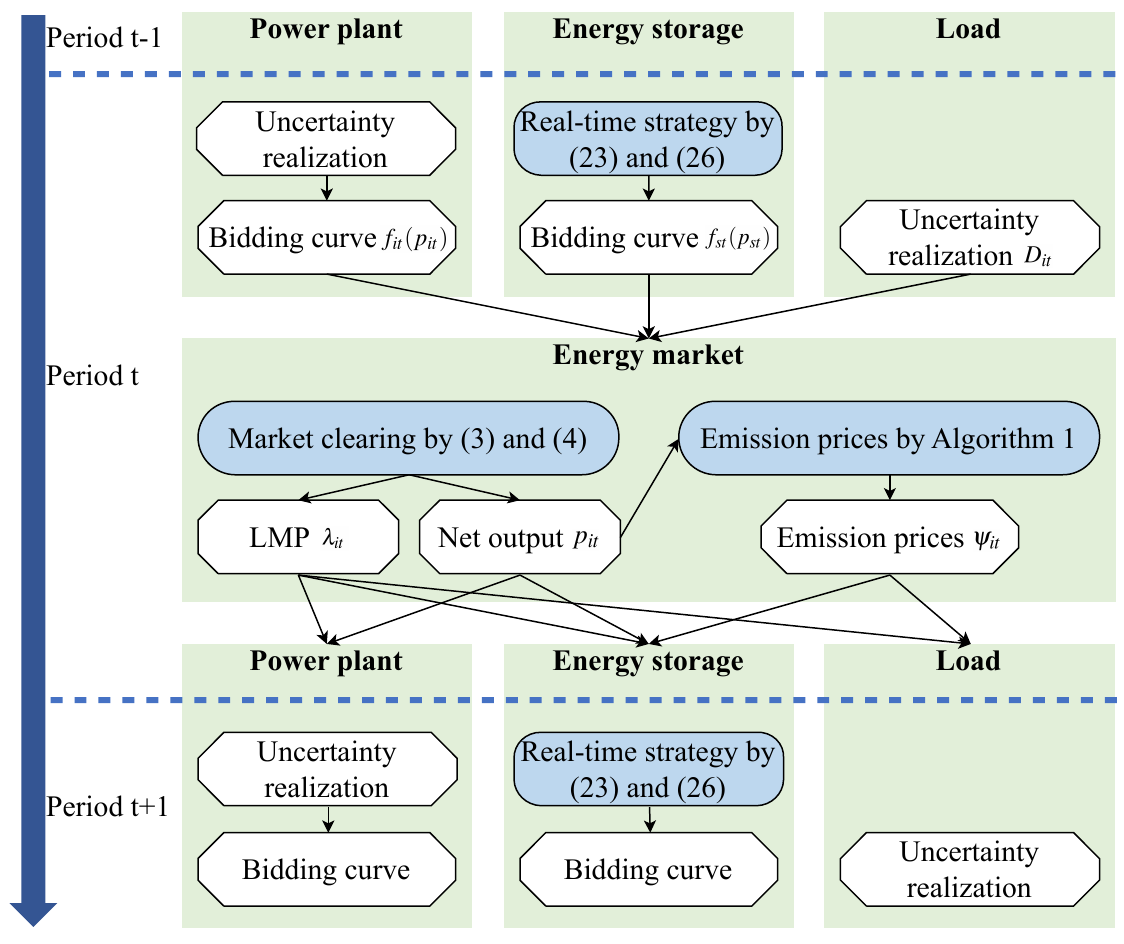}
\caption{The overall procedure of the proposed electricity market.}
\label{fig:framework}
\end{figure}

The overall procedure of the proposed electricity market is illustrated in Fig. \ref{fig:framework}. At the beginning of period $t$, uncertain renewable power plants and loads observe their actual power output bounds/demands. Power plant $i \in S_G$ submits the bidding curve $f_{it}(p_{it})$ in \eqref{eq:bidding-plant} to the market operator. Based on the emission price $\psi_{s(t-1)}$ and the virtual queue $q_{st}$, ES $s \in S_S$ submits the bidding curve in \eqref{eq:bidding-ES} and the bounds given by \eqref{eq:bounds-ES} to the market operator.
Load $i \in S_B$ reports its demand power $D_{it}$. Then the market operator solves the OPF problem \eqref{eq:market-LP} to obtain the net output $p_{it}, i \in S_G \cup S_S$ for power plants and ESs. In addition, the LMPs $\lambda_{it}, i \in S_B$ are calculated by \eqref{eq:LMP}. Subsequently, emission prices $\psi_{it}, i \in S_B$ are calculated using Algorithm \ref{alg:allocation} and ESs and loads pay for half of the total emissions. After period $t$ ends, period $t+1$ starts, and the above process repeats.

\emph{Remark:} The strategies of ES in the day-ahead market \cite{wang2017look} and the ancillary service market \cite{cho2015energy} are out of the scope of this paper. 
Nonetheless, the proposed ES bidding strategy can be adapted to consider the previous clearing results in the day-ahead energy market and the ancillary service market. Based on the previous clearing results, the equivalent adjustable region of ES can be obtained, which contains the power output range that the ES can provide in the real-time energy market. This range information can be added to the power output bounds of the bidding curve.

\section{Case Studies}
\label{sec-5}

In this section, we first test the proposed method using a modified IEEE 30-bus case. To demonstrate its effectiveness and advantages, the proposed method is compared with multiple existing methods. The impact of different factors is also analyzed. The scalability is examined by a modified IEEE 118-bus case. All the experiments are done on a laptop with an Intel i7-12700H processor and 16 GM RAM. Linear programs are solved by Gurobi 9.5.

\subsection{Performance Evaluation}

\begin{figure}[t]
\centering
\includegraphics[width=1.0\columnwidth]{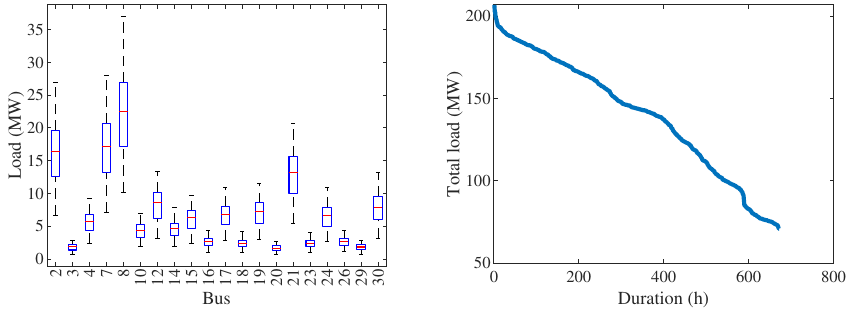}
\caption{The boxplot of nodal load demand (left) and the load duration curve of the total load (right) in the modified IEEE 30-bus system.}
\label{fig:load}
\end{figure}

As a benchmark, the modified IEEE 30-bus case is tested. 
There are 6 fossil fuel generators, whose parameters are in TABLE \ref{table:generator}. The load data are shown in Fig. \ref{fig:load}. 
A $100$-MW PV station and a $100$-MW wind power plant are connected to buses 6 and 15, respectively. Two ESs are equipped at buses 15 and 18, respectively, whose parameters are listed in TABLE \ref{table:parameter}. Other data can be found in \cite{xie2023github}.
We run the simulation over $28$ days divided into $672$ periods (1 h each). It takes about 400 s to output the result. To evaluate its performance, comparisons are conducted in the following.


\begin{table}[!ht]
\scriptsize
\renewcommand{\arraystretch}{1.3}
\caption{Parameters of fossil fuel generators in the modified IEEE 30-bus system}
\label{table:generator}
\centering
\begin{tabular}{ll}
\hline
Parameter & Value \\
\hline
Bus & (1, 2, 22, 27, 23, 13) \\
Fuel cost coefficient (\$/kWh) & (0.047, 0.055, 0.055, 0.047, 0.055, 0.047) \\
Emission coefficient (kgCO$_2$/kWh) & (0.9, 0.8, 0.8, 0.2, 0.3, 0.3) \\
Maximum output (MW) & (80, 80, 50, 55, 30, 40) \\
\hline
\end{tabular}
\end{table}

\begin{table}[!t]
\scriptsize
\renewcommand{\arraystretch}{1.3}
\caption{Other parameters in the modified IEEE 30-bus system}
\label{table:parameter}
\centering
\begin{tabular}{cc|cc}
\hline
Parameter & Value & Parameter & Value \\
\hline
$T$ & 672 & $\tau$ & $1$ h \\
$\epsilon$ & $0.0001$ \$/kgCO$_2$ & $\kappa$ & $0.05$ \$/kgCO$_2$ \\
$\delta$ & $0.002$ & $(N_{15},N_{18})$ & $(50,50)$ \\
$(\eta^c,\eta^d)$ & $(0.95,0.95)$ & $(P_{15}^{max},P_{18}^{max})$ & $(4,4)$ MW \\
$(\underline{E}_{15},\underline{E}_{18})$ & $(4,2)$ MWh & $(\overline{E}_{15},\overline{E}_{18})$ & $(36,18)$ MWh \\
\hline
\end{tabular}
\end{table}


\subsubsection{Cost and Emission Comparison}

We compare the generation cost and carbon emission of four markets with/without ESs and carbon emission allocation. The settings and results are given in TABLE \ref{table:comparison}. 
The renewable curtailment rate is calculated using the optimal solution of the market clearing OPF problem \eqref{eq:market-LP}, which equals:
\begin{align}
    \frac{\sum_{t = 1}^T \sum_{i \in S_R} (\overline{P}_{it}-p_{it})}{\sum_{t = 1}^T \sum_{i \in S_R} \overline{P}_{it}}, \nonumber
\end{align}
where $S_R$ is the index set of renewable power plants, $\overline{P}_{it}$ is maximum renewable power, and $p_{it}$ is the power output in the market clearing result, so $\overline{P}_{it}-p_{it}$ is the curtailment power of power plant $i$ in period $t$.

The only difference between the settings of case Proposed and case A1 in TABLE \ref{table:comparison} is that case A1 does not employ carbon emission allocation while case Proposed does. The total emission in case Proposed is 43.1\% lower than that in case A1. The main reason is that a large amount of power generation is shifted from carbon-intensive but low-price power plants (such as coal-fired units) to greener power plants with a relatively high unit generation cost (such as gas-fired units). The carbon emission allocation mechanism changes the market dynamics by adding emission costs and making carbon-intensive power plants less competitive. Therefore, compared with case A1, the total generation cost increases and the total emission decreases in case Proposed.
The benefit of ES participation can be observed by comparing Proposed and A2, where the total generation cost, total emission, and renewable curtailment are reduced by 1.63\%, 1.66\%, and 43.4\%, respectively.

\begin{table}[!ht]
\scriptsize
\renewcommand{\arraystretch}{1.3}
\caption{Results with/without ESs and carbon emission allocation}
\label{table:comparison}
\centering
\begin{tabular}{ccccc}
\hline
Case & Proposed & A1 & A2 & A3 \\
\hline
ESs & \checkmark & \checkmark & $\times$ & $\times$ \\
Carbon emission allocation & \checkmark & $\times$ & \checkmark & $\times$ \\
Total generation cost (\$/h) & $3387$ & $3121$ & $3443$ & $3173$ \\
Total emission (kgCO$_2$/h) & $30546$ & $53701$ & $31063$ & $54457$ \\
Renewable curtailment & $1.84$\% & $1.84$\% & $3.25$\% & $3.25$\% \\
\hline
\end{tabular}
\end{table}

\begin{figure}[t]
\centering
\includegraphics[width=1.0\columnwidth]{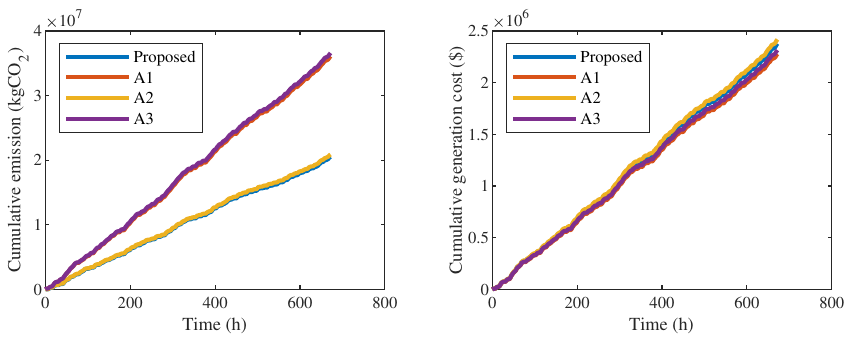}
\caption{Cumulative emission (left) and cumulative generation cost (right) in different cases.}
\label{fig:cumulative-operation}
\end{figure}

\begin{figure}[t]
\centering
\includegraphics[width=0.85\columnwidth]{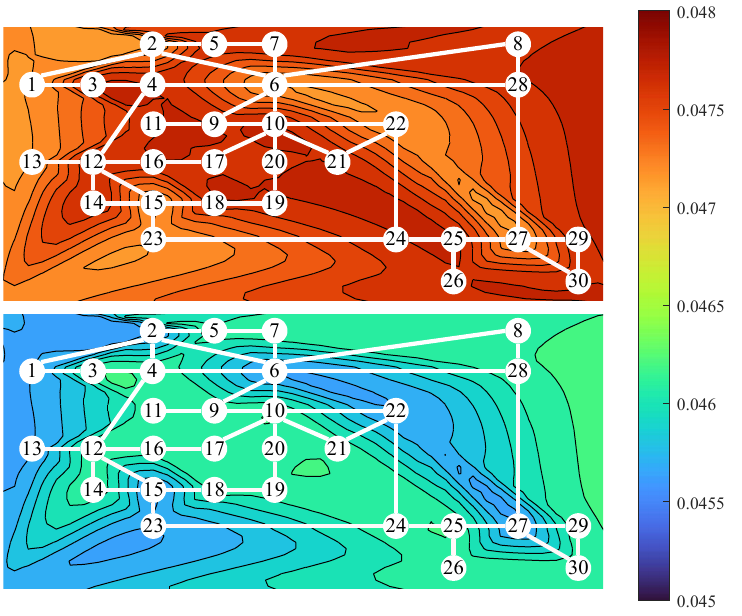}
\caption{Electricity prices (\$/kWh) in the modified IEEE 30-bus case in case Proposed (the first contour map) and case A1 (the second contour map).}
\label{fig:electricity-price}
\end{figure}

The cumulative emission and generation cost curves are depicted in Fig. \ref{fig:cumulative-operation}. The carbon emission allocation significantly decreases the total emission, whereas the ES participation slightly decreases it. The carbon emission allocation increases the total generation cost by about 9\%, and the ES participation decreases it by about 2\%. The nodal electricity prices in case Proposed and case A1 in period 50 are compared in Fig. \ref{fig:electricity-price}, which shows that they have similar distribution patterns, but the electricity prices are higher in case Proposed. This is because when fossil fuel generators have to pay for their emission costs, their bidding cost curves become higher and then the electricity prices increase.

The approximation error of the emission prices depends on the specific case. If the emission prices change slowly across time or the emission cost coefficient $\kappa$ is small, the revenue error introduced by approximating the emission prices with the prices from the previous period will be relatively small. For example, in the benchmark case, this revenue error does not exceed 1\%. However, it can be larger in other cases. Besides directly approximating emission prices by the previous values, the ES operator can utilize historical price data and employ price forecasting methods for better approximation performance.

\subsubsection{Energy Storage Revenue Comparison}
To show the effectiveness of the proposed real-time ES bidding strategy, we compare it with three alternatives as follows. 
The impacts of the ES operation strategy on the market price are not considered in the three traditional methods, and we focus on the ES operation strategy, i.e., the power-price function, rather than the bidding strategy in this experiment. Therefore, we fix the prices to those in the benchmark case and assume that the price curve will not be affected by the ES operation strategy. In the offline method, all the prices are assumed to be known, including future prices. 
\begin{itemize}
    \item B1: \emph{Real-time Strategy Based on Traditional Lyapunov Optimization}. Different from the proposed method, B1 minimizes an upper bound of the drift-plus-penalty $\Delta_{st}$ rather than the exact one in \eqref{eq:drift-plus-penalty}, i.e., the objective of \textbf{P2} is replaced by
    
\vspace{-1em}

{\begin{small}
\begin{equation}
    \min_{p_{st}^c,p_{st}^d}~ C_0 + (p_{st}^c \tau \eta_s^c - p_{st}^d \tau / \eta_s^d)q_{st}  + V_s(-\gamma_{st} (p_{st}^d-p_{st}^c) \tau), \nonumber
\end{equation}
\end{small}}

where $C_0$ is a constant 
and the parameters are set as
\begin{align}
    & V_s = \frac{\overline{E}_s - \underline{E}_s - P_s^{max} \tau \eta_s^c - P_s^{max} \tau / \eta_s^d}{\overline{\gamma}_s \eta_s^d - \underline{\gamma}_s/\eta_s^c}, \nonumber \\
    & E_s = \overline{E}_s + V_s \underline{\gamma}_s / \eta_s^c - P_s^{max} \tau \eta_s^c. \nonumber
\end{align}
\item B2: \emph{Simple Strategy}. A simple strategy with lower and upper price thresholds is used. The ES charges/discharges with the maximum feasible power if the combined price is below $0.02$ \$/kWh/above $0.05$ \$/kWh. 
\item B3: \emph{Strategy By the Offline Model in \eqref{eq:offline}}.
\end{itemize}
 
\begin{figure}[t]
\centering
\includegraphics[width=0.70\columnwidth]{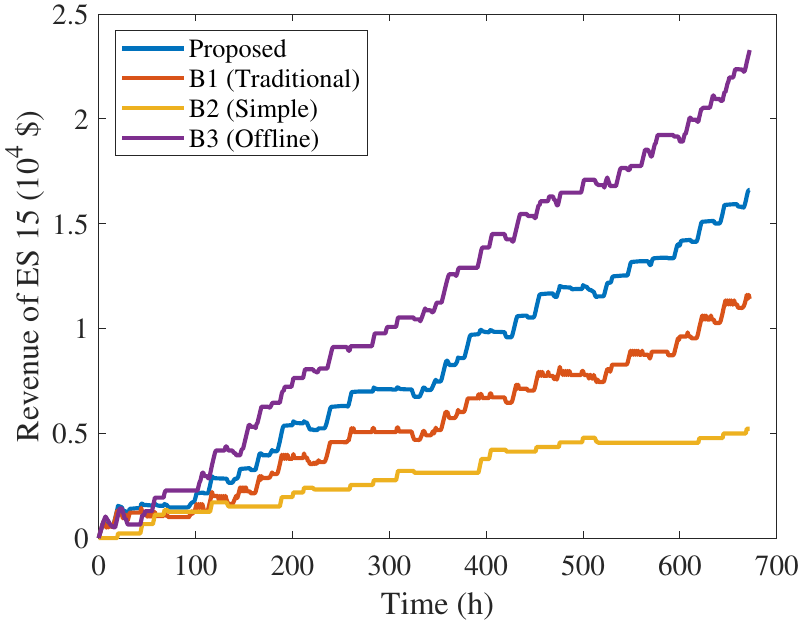}
\caption{Revenue of ES 15 by different ES operation methods.}
\label{fig:revenue}
\end{figure}

We use the ES at bus 15 (ES 15 for short) as an example.
Its revenue curves under case Proposed and B1-B3 are compared in Fig. \ref{fig:revenue}. When the ES charges, it pays the bill; when it discharges, it gains income. Therefore, the revenue curves oscillate up and down. 
The least squares technique is applied to fit the revenue curves by linear functions, whose slopes can be regarded as approximate revenue rates. The revenue rates of the proposed method and B1-B3 are $23.50$ \$/h, $15.65$ \$/h, $7.48$ \$/h, and $33.14$ \$/h, respectively. The offline method assumes complete knowledge of future prices to find the optimal strategy, so it has the highest revenue rate. However, it is not practical. 
The proposed method and B1 perform much better than B2, and achieve about $70.9$\% and $47.2$\% of the offline revenue rate. This shows that the proposed method is notably more effective than the traditional Lyapunov optimization-based method (B1).

\subsubsection{Comparison of Carbon Emission Allocation Methods}

To demonstrate the effectiveness of the proposed carbon emission allocation mechanism based on Aumann-Shapley prices, we compare it with the CEF method \cite{kang2015carbon}.

\begin{figure}[t]
\centering
\includegraphics[width=0.85\columnwidth]{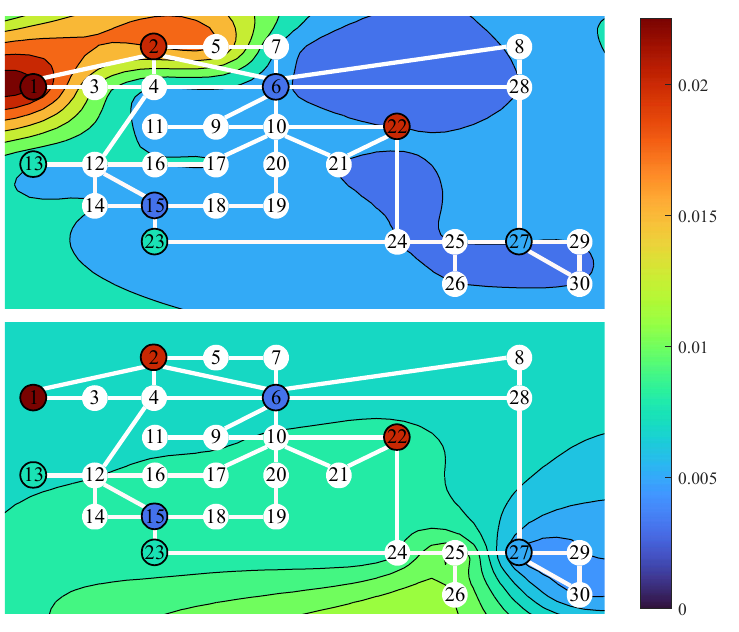}
\caption{Emission prices (\$/kWh) in the modified IEEE 30-bus case by the CEF method (the first contour map) and the proposed method (the second contour map). The buses connecting generators are also painted to show their emission prices (i.e., $\kappa \Psi_i/2$ for power plant $i$).}
\label{fig:CEFintensity}
\end{figure}

We first investigate the emission prices deduced by the two methods. To be fair, we test the two methods using a system without ES so that the total emissions within each period under the two methods are equal. 
The emission price at bus $i$ by the CEF method is half of the CEF intensity flowing out of the bus. The emission prices in period 50 at different buses by the two methods are plotted using contour maps in Fig. \ref{fig:CEFintensity}. The prices vary more significantly if the CEF method is used. Bus 26 has the highest emission price under the proposed method, while the CEF method gives a relatively low emission price for bus 26. To see which method can better reflect the contribution of load demand at bus 26 to carbon emissions, we test the sensitivity of the total emission toward the demand at bus 26. We can get that $(\partial \mathcal{E}/\partial D_{26})(D^*)/\kappa \approx 1.85$ kgCO$_2$/kWh, which is the highest among all the buses. This shows the proposed method is more effective.

\begin{figure}[t]
\centering
\includegraphics[width=1.0\columnwidth]{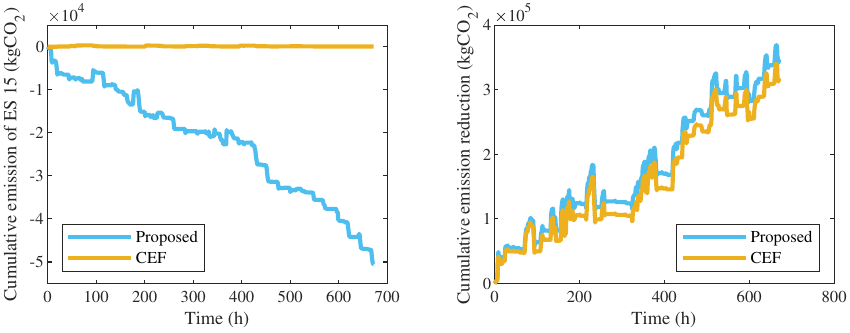}
\caption{Cumulative allocated emission of ES 15 (left) and cumulative total emission reduction (right) by different methods.}
\label{fig:emission}
\end{figure}

To further test the performance of the two methods in reducing total emissions, we test the two methods using a system with ESs. The CEF method considering ES in \cite{wang2021optimal} and \cite{yang2023improved} is applied. The cumulative allocated emissions of ES 15 are compared in Fig. \ref{fig:emission} (left). Since ES 15 helps with the system emission reduction, its cumulative emission decreases (with an emission rate of $-62.68$ kgCO$_2$/h) in the proposed method. For the CEF method, the cumulative allocated emission is slowly increasing (with an emission rate of $0.12$ kgCO$_2$/h) due to the energy loss in the charging and discharging processes. The system total emission is compared with the no ES case (A2 in TABLE \ref{table:comparison}), and the cumulative reduction values are shown in Fig. \ref{fig:emission} (right). The system total emission by the CEF method is $30590$ kgCO$_2$/h and the system total emission reduction is $473$ kgCO$_2$/h. In contrast, the reduction by the proposed method is $517$ kgCO$_2$/h, which is about $9.3$\% higher than the CEF method. This is because the proposed method directly measures the impact of ES power on the total emission, while the CEF method fails to reflect the impact of ES discharging on carbon intensity as it depends only on the inflow.

\subsubsection{Accuracy Comparison}
To show the accuracy of Algorithm \ref{alg:allocation} for calculating the Aumann-Shapley prices-based carbon emission allocation, we compare it with the following two numerical calculation methods:
\begin{itemize}
\item C1: The partial derivatives and the integrals in \eqref{eq:emission-allocation} are calculated using numerical estimations \cite{chen2018method,nan2022bi,nan2022hierarchical,zhou2019cooperative}.
\item C2: The partial derivatives in \eqref{eq:emission-allocation} are calculated using the analytical expression in \eqref{eq:partial-emission}, while the integrals are computed numerically.
\end{itemize}
The results in period 50 are listed in TABLE \ref{table:numerical}. The number of sample points on the segment from $0$ to $\tilde{D}^*$ is in the second column. In the proposed method, the sample number equals the number of iterations of Algorithm \ref{alg:allocation}. The accuracy is measured by the cost-sharing error, which is defined as the relative error from the sum of the allocated emissions to the total emission being allocated. The results show that method C1 has the largest errors and the longest computation time. Using the analytical expression of partial derivatives, method C2 has a better performance than method C1 but still needs much more sample points and a longer computation time than the proposed method to achieve satisfactory accuracy. Therefore, the proposed algorithm is more precise.
\begin{table}[!t]
\scriptsize
\renewcommand{\arraystretch}{1.3}
\caption{Comparison of different carbon emission allocation calculation methods}
\label{table:numerical}
\centering
\begin{tabular}{cccc}
\hline
Method & Sample number & Cost-sharing error & Computation Time (s) \\
\hline
C1 & $100$ & $4.64$\% & $159$ \\
C1 & $1000$ & $3.19$\% & $1680$ \\
C2 & $100$ & $2.74$\% & $6.91$ \\
C2 & $1000$ & $0.02$\% & $107$ \\
Proposed & $4$ & $0.00$\% & $0.37$ \\
\hline
\end{tabular}
\end{table}

\subsection{Impact of Some Factors}

\begin{figure}[t]
\centering
\includegraphics[width=1.0\columnwidth]{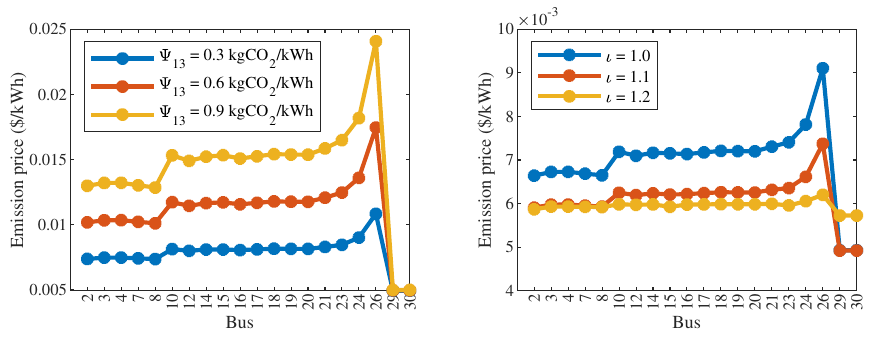}
\caption{Emission prices (\$/kWh) in the modified IEEE 30-bus case under different fossil fuel generator emission coefficient $\Psi_{13}$ (left) and under different multiple $\iota$ of transmission line capacity (right).}
\label{fig:sensitivity}
\end{figure}

We test and analyze the impacts of the main parameters. First, we change the emission coefficient $\Psi_{13}$ of the fossil fuel generator at bus 13. The emission prices in period 50 are shown on the left side in Fig. \ref{fig:sensitivity}. Only the buses with nonzero demand are drawn. As $\Psi_{13}$ increases, the emission prices of most buses increase. The emission prices at bus 29 and bus 30 do not change much, because almost all the electricity flowing into the two buses comes from the generator at bus 27.

Then we investigate the impacts of the transmission line capacity. We multiply the transmission line capacity $F_l,\forall l \in S_L$ by a constant $\iota$. The emission prices under different $\iota$ values are shown on the right side in Fig. \ref{fig:sensitivity}. As the transmission line capacity becomes larger, the network congestion is alleviated. Thus, the values of the partial derivatives $\partial \mathcal{E}/\partial D_i, i \in S_B$ have a convergence tendency, which indicates that the demand location is making a smaller difference in the emission prices. This is exactly what is found in Fig. \ref{fig:sensitivity}.

We then investigate the Lyapunov optimization parameter $V_s$. 
Take ES 15 as an example, we change its parameter $V_{15}$ to 0.1, 0.4, 0.7, and 1.0 times the proposed value in \eqref{eq:V-E-value} and the corresponding $E_{15}$ is set as the average value of the lower and upper bounds in \eqref{eq:V-E-range-2}.
The results are shown in TABLE \ref{table:Vs}. A larger $V_s$ leads to a higher revenue rate and a lower system total emission, which is consistent with the theoretical analysis in Theorem \ref{thm:performance}.
\begin{table}[!t]
\scriptsize
\renewcommand{\arraystretch}{1.3}
\caption{Results under different Lyapunov optimization parameter $V_{15}$}
\label{table:Vs}
\centering
\begin{tabular}{ccccc}
\hline
Multiple of $V_{15}$ & $0.1$ & $0.4$ & $0.7$ & $1.0$ \\
\hline
Revenue of ES 15 (\$/h) & $2.53$ & $9.49$ & $17.87$ & $23.50$ \\
Emission of ES 15 (kgCO$_2$/h) & $-2.47$ & $-16.46$ & $-44.84$ & $-62.68$ \\
System emission (kgCO$_2$/h) & $30845$ & $30793$ & $30675$ & $30546$ \\
\hline
\end{tabular}
\end{table}

The emission cost coefficient $\kappa$ is also tested and the results are in TABLE \ref{table:kappa}. A larger $\kappa$ means the emission weighs more in the total cost, so the system emission decreases as $\kappa$ increases. 
An interesting finding is that when $\kappa$ increases, the ES revenue rate may increase because ES has the opportunity to earn a higher profit from contributing to carbon emission reduction. The change of $\kappa$ has impacts on the nodal emission price $\psi_{st}$, the bidding cost curve $f_{it} (p_{it})$, and the ES emission cost $-\psi_{st} p_{st}$. When $\kappa$ increases, the ES has a better opportunity to earn from contributing to carbon emission reduction, but its revenue may decrease if it relies mainly on the electricity price to make profits and the changing trends of electricity prices and emission prices are different.
\begin{table}[!t]
\scriptsize
\renewcommand{\arraystretch}{1.3}
\caption{Results under different cost coefficient $\kappa$ of emission}
\label{table:kappa}
\centering
\begin{tabular}{ccccc}
\hline
$\kappa$ (\$/kgCO$_2$) & $0$ & $0.02$ & $0.05$ & $0.10$ \\
\hline
Revenue of ES 15 (\$/h) & $21.53$ & $22.34$ & $23.50$ & $26.05$ \\
Revenue of ES 18 (\$/h) & $9.19$ & $9.68$ & $10.31$ & $10.98$ \\
System emission (kgCO$_2$/h) & $53701$ & $53530$ & $30546$ & $28482$ \\
\hline
\end{tabular}
\end{table}

\subsection{Scalability}

The proposed method is further tested on a modified IEEE 118-bus case to demonstrate its scalability, where the data are in \cite{xie2023github}. There are $56$ power plants and $99$ loads in the system. The average computation time and iterations of Algorithm \ref{alg:allocation} of one period under different numbers of ESs are listed in TABLE \ref{table:time}. It shows that the proposed method is computationally efficient enough for a real-time electricity market. 
\begin{table}[!ht]
\scriptsize
\renewcommand{\arraystretch}{1.3}
\caption{Average computation time/average iterations of Algorithm \ref{alg:allocation} under different settings}
\label{table:time}
\centering
\begin{tabular}{cccc}
\hline
Number of ESs & $2$ & $8$ & $16$ \\
\hline
IEEE 30-bus & $0.65$ s/$4.25$ & $0.68$ s/$4.09$ & $0.77$ s/$4.09$ \\
IEEE 118-bus & $2.23$ s/$10.83$ & $2.69$ s/$10.83$ & $3.60$ s/$10.82$ \\
\hline
\end{tabular}
\end{table}

\section{Conclusion}
\label{sec-6}

This paper first proposes an electricity market with an Aumann-Shapley prices-based carbon emission allocation mechanism. A parametric linear programming-based method is proposed to calculate the carbon emission allocation. Then, a real-time (online) bidding strategy for ES to participate in the market is developed based on Lyapunov optimization.
The main findings of the case studies include:
\begin{itemize}
    \item The proposed Aumann-Shapley prices-based emission prices can better encourage ESs to help with system emission reduction than the traditional CEF method.
    \item The proposed emission calculation method is more accurate and efficient than the existing numerical methods.
    \item The proposed real-time ES strategy using exact drift-plus-penalty minimization leads to higher revenue rates than applying the traditional Lyapunov optimization technique.
\end{itemize}

  
\ifCLASSOPTIONcaptionsoff
\newpage
\fi

\appendices
\makeatletter
\@addtoreset{equation}{section}
\@addtoreset{theorem}{section}
\makeatother

\setcounter{equation}{0}  
\renewcommand{\theequation}{A.\arabic{equation}}
\renewcommand{\thetheorem}{A.\arabic{theorem}}
\section{Properties of the Proposed Carbon Emission Allocation Mechanism}
\label{appendix-A}
\subsection{Proof of Proposition \ref{prop:cost-sharing}}

Let $\theta \triangleq (p, d)$. Then the standard compact form of the OPF problem \eqref{eq:market-fixed} is
\begin{align}
    \min_x~ & C^\top x \nonumber \\
    \label{eq:market-fixed-compact}
    \mbox{s.t.}~ & Ax = Q \theta + H, x \geq 0 
\end{align}
which can be regarded as a multiparametric linear programming problem with the multi-dimensional parameter $\theta$, where $\theta$ only appears in the right-hand-side coefficient of the constraint. In addition, the emission $\mathcal{E} (\theta) = K^\top x$ for some coefficient $K$.

According to the multiparametric linear programming theory \cite{gal1972multiparametric}, the feasible range of $\theta$ (within which \eqref{eq:market-fixed-compact} has a feasible solution) is a convex set. Problem \eqref{eq:market-fixed-compact} is feasible for $\theta = 0$ and $\theta = \Theta^* \triangleq (P^*,D^*)$. Thus, it is feasible for $\theta = y\Theta^*, \forall y \in [0,1]$. By the physical meaning, variable $x$ in \eqref{eq:market-fixed-compact} is bounded, so for $y\Theta^*, \forall y \in [0,1]$, the emission $\mathcal{E}(y\Theta^*)$ has a finite value.

Again by multiparametric linear programming \cite{gal1972multiparametric}, the feasible range of $\theta$ can be divided into a finite number of critical regions. In a critical region, an optimal basis remains the same, so the optimal solution and the emission cost $\mathcal{E}(\theta)$ is affine in $\theta$. In addition, every critical region is a closed polyhedron. The straight line from $\theta = 0$ to $\theta = \Theta^*$ can be divided into a finite number of segments in different critical regions. Suppose $0 = y_0 < y_1 < \dots < y_M = 1$ and for any $1 \le m \le M$, the segment from $y_{m-1} \Theta^*$ to $y_m \Theta^*$ is in the same critical region. 

For any $m$ with $1 \leq m \leq M$, assume the optimal basis in the critical region $\Omega_{m-1}$ containing the segment from $y_{m-1}\Theta^*$ to $y_m\Theta^*$ is $A_{B_{m-1}}$, then $x = (x_B, x_N)$ with $x_B = A_{B_{m-1}}^{-1} (Q \theta + H)$ and $x_N = 0$ is an optimal solution for $\theta$ in this critical region. Thus, $\mathcal{E}(\theta) = K_{B_{m-1}}^\top A_{B_{m-1}}^{-1} (Q \theta + H), \theta \in \Omega_{m-1}$, which is an affine function and hence smooth on $\Omega_{m-1}$. Therefore, by Newton-Leibniz theorem and the chain rule of derivatives, 
\begin{align}
    \mathcal{E}(y_m \Theta^*) - \mathcal{E}(y_{m-1} \Theta^*) & = \int_{y_{m-1}}^{y_m} \frac{d\mathcal{E}(y\Theta^*)}{dy} dy \nonumber \\
    & = \int_{y_{m-1}}^{y_m} \left( \sum_j \Theta_j^* \cdot \frac{\partial \mathcal{E}}{\partial \theta_j} (y\Theta^*) \right) dy \nonumber \\
    & = \sum_j \Theta_j^* \int_{y_{m-1}}^{y_m} \frac{\partial \mathcal{E}}{\partial \theta_j} (y\Theta^*) dy. \nonumber
\end{align}
Then
\begin{align}
    \mathcal{E}(\Theta^*) - \mathcal{E}(0) & = \sum_{m = 1}^M (\mathcal{E}(y_m \Theta^*) - \mathcal{E}(y_{m-1} \Theta^*)) \nonumber \\
    & = \sum_{m = 1}^M \sum_j \Theta_j^* \int_{y_{m-1}}^{y_m} \frac{\partial \mathcal{E}}{\partial \theta_j} (y\Theta^*) dy \nonumber \\
    & = \sum_j \Theta_j^* \int_0^1 \frac{\partial \mathcal{E}}{\partial \theta_j} (y\Theta^*) dy \nonumber \\
    & = \sum_j \int_0^{\Theta_j^*} \frac{\partial \mathcal{E}}{\partial \theta_j} (\frac{y}{\Theta_j^*}\Theta^*) dy \nonumber\\
    & = \sum_{s \in S_S} \mathcal{E}_s(\Theta^*) + \sum_{i \in S_B} \mathcal{E}_i(\Theta^*). \nonumber
\end{align}
This completes the proof.

\subsection{Other Properties}
Apart from the cost-sharing property, the proposed carbon emission allocation mechanism also possesses other nice properties listed below, which are direct consequences of the definitions in \eqref{eq:emission-allocation} and \eqref{eq:emission-price}. For conciseness, they are stated in terms of loads but also apply to ESs.
\begin{itemize}
    \item Scale invariance: The allocation results are independent of the units.
    \item Monotonicity: If $\partial\mathcal{E}/\partial D_i(D) \geq \partial\mathcal{E}/\partial D_j(D), \forall D$, then $\psi_i(D) \geq \psi_j(D), \forall D$ holds. In other words, the load that always has a larger influence on the total emission will receive a higher emission price.
    \item Additivity: If $\mathcal{E}(D) = \tilde{\mathcal{E}}(D) + \check{\mathcal{E}}(D), \forall D$, then the allocation results satisfy $\mathcal{E}_i(D) = \tilde{\mathcal{E}}_i(D) + \check{\mathcal{E}}_i(D), \forall D, \forall i$. According to additivity, the allocation results will remain the same as that in \eqref{eq:emission-allocation} if we allocate the emission of each power plant among loads and then sum them up.
    \item Consistency: If $\mathcal{E}(D) = \tilde{\mathcal{E}}(\sum_{i \in \mathcal{I}} D_i, D_{i'}, i' \notin \mathcal{I}), \forall D$, then $\psi_j(D) = \tilde{\psi}_{\mathcal{I}}(\sum_{i \in \mathcal{I}} D_i, D_{i'}, i' \notin \mathcal{I}), \forall D, \forall j \in \mathcal{I}$. The meaning is that if the emission to be allocated is a function of the sum of some loads, then these loads can be merged before the allocation.
\end{itemize}

\setcounter{equation}{0}  
\renewcommand{\theequation}{B.\arabic{equation}}
\renewcommand{\thetheorem}{B.\arabic{theorem}}
\section{Proof of Theorem \ref{thm:feasible}}
\label{appendix-B}

The constraint $\underline{E}_s \leq e_{st} \leq \overline{E}_s$ is equivalent to 
\begin{align}
\label{eq:q-constraint}
    \underline{E}_s - E_s \leq q_{st} \leq \overline{E}_s - E_s.
\end{align}
We prove the conclusion by mathematical induction: Assume \eqref{eq:q-constraint} holds and prove $\underline{E}_s - E_s \leq q_{s(t+1)} \leq \overline{E}_s - E_s$.

Combine $q_{s(t+1)} = q_{st}+p_{st}^c \tau \eta_s^c - p_{st}^d \tau / \eta_s^d$ and $p_{st} = p_{st}^d - p_{st}^c$ with \eqref{eq:strategy-gamma}, we have
\begin{footnotesize}
\begin{align}
    & q_{s(t+1)} =  \nonumber \\
    & \left\{
    \begin{array}{ll}
        q_{st}+P_s^{max} \tau\eta_s^c , &\text{if}~ q_{st} \in (-\infty, - \frac{V_s \gamma_{st}}{\eta_s^c} - P_s^{max} \tau\eta_s^c ], \\
        - \frac{V_s \gamma_{st}}{\eta_s^c}, &\text{if}~ q_{st} \in [- \frac{V_s \gamma_{st}}{\eta_s^c} - P_s^{max} \tau \eta_s^c, - \frac{V_s \gamma_{st}}{\eta_s^c}], \\
        q_{st}, &\text{if}~ q_{st} \in [- \frac{V_s \gamma_{st}}{\eta_s^c}, - V_s \gamma_{st} \eta_s^d], \\
        -V_s \gamma_{st} \eta_s^d, &\text{if}~ q_{st} \in [- V_s \gamma_{st} \eta_s^d, - V_s \gamma_{st} \eta_s^d + \frac{P_s^{max} \tau}{\eta_s^d}], \\
        q_{st}-P_s^{max}\tau/\eta_s^d, &\text{if}~ q_{st} \in [ - V_s \gamma_{st} \eta_s^d + \frac{P_s^{max} \tau}{\eta_s^d}, +\infty).
    \end{array}
    \right. \nonumber
\end{align}
\end{footnotesize}
Therefore,
\begin{small}
\begin{align}
& q_{s(t+1)} \in \nonumber \\
     & \left\{
    \begin{array}{ll}
        [q_{st}, - \frac{V_s \gamma_{st}}{\eta_s^c}), &\text{if}~ q_{st} \in (-\infty, - \frac{V_s \gamma_{st}}{\eta_s^c} - P_s^{max} \tau \eta_s^c), \\
        \{- \frac{V_s \gamma_{st}}{\eta_s^c}\}, &\text{if}~ q_{st} \in [- \frac{V_s \gamma_{st}}{\eta_s^c} - P_s^{max} \tau \eta_s^c, - \frac{V_s \gamma_{st}}{\eta_s^c}], \\
        \{q_{st}\}, &\text{if}~ q_{st} \in [- \frac{V_s \gamma_{st}}{\eta_s^c}, - V_s \gamma_{st} \eta_s^d], \\
        \{-V_s \gamma_{st} \eta_s^d\}, &\text{if}~ q_{st} \in [- V_s \gamma_{st} \eta_s^d, - V_s \gamma_{st} \eta_s^d + \frac{P_s^{max} \tau}{\eta_s^d}], \\
        (-V_s \gamma_{st} \eta_s^d, q_{st}], &\text{if}~ q_{st} \in ( - V_s \gamma_{st} \eta_s^d + \frac{P_s^{max} \tau}{\eta_s^d}, +\infty).
    \end{array}
    \right. \nonumber
\end{align}
\end{small}
By \eqref{eq:q-constraint}, we only need to prove that $\forall \gamma_{st} \in [\underline{\gamma}_s, \overline{\gamma}_s],$
\begin{align}
    \label{eq:gamma-constraint-1}
    - V_s \gamma_{st}/\eta_s^c \leq \overline{E}_s - E_s,~ - V_s \gamma_{st} \eta_s^d \geq \underline{E}_s - E_s.
\end{align}
Recall that $V_s > 0$ and $0 \leq \underline{\gamma}_{st} < \eta_c \eta_d \overline{\gamma}_{st}$, then we have
\begin{subequations}
\begin{align}
    \eqref{eq:gamma-constraint-1} & \iff - V_s \underline{\gamma}_s/\eta_s^c \leq \overline{E}_s - E_s,~ - V_s \overline{\gamma}_s \eta_s^d \geq \underline{E}_s - E_s \\
    \label{eq:gamma-constraint-2}
    & \iff \underline{E}_s + V_s \overline{\gamma}_s \eta_s^d \leq E_s \leq \overline{E}_s + V_s \underline{\gamma}_s / \eta_s^c.
\end{align}
\end{subequations}

When \eqref{eq:V-E-range-1} holds, $\underline{E}_s + V_s \overline{\gamma}_s \eta_s^d \leq \overline{E}_s + V_s \underline{\gamma}_s / \eta_s^c$, so there exists $E_s$ satisfying \eqref{eq:gamma-constraint-2}, which is the same as \eqref{eq:V-E-range-2}. Therefore, $\underline{E}_s - E_s \leq q_{s(t+1)} \leq \overline{E}_s - E_s$ and the conclusion follows from mathematical induction.

\setcounter{equation}{0}  
\renewcommand{\theequation}{C.\arabic{equation}}
\renewcommand{\thetheorem}{C.\arabic{theorem}}
\section{Proof of Theorem \ref{thm:performance}}
\label{appendix-C}

Denote the strategy in \eqref{eq:strategy-gamma}, the corresponding virtual queue, and Lyapunov drift by $p_{st}^*$, $q_{st}^*$, and $\Delta_{st}^*$ for any period $t$, respectively. According to Theorem \ref{thm:feasible}, this strategy is feasible for problem \eqref{eq:v_0}. Then by the optimality of $v_0^*$ in problem \eqref{eq:v_0}, we have $v^* \leq v_0^*$. Note that $v_1^* \geq v_0^*$ because problem \eqref{eq:v_1} is obtained from problem \eqref{eq:v_0} by relaxing a constraint, so we only need to show that
\begin{align}
\label{eq:v_1-inequality}
    v^* \geq v_1^*-(P_s^{max} \tau)^2 / (2V_s(\eta_s^d)^2).
\end{align}

According to Lyapunov optimization theory \cite{neely2010stochastic}, for any $\tilde{\delta} > 0$, there is a so-called $\omega$-only policy $\tilde{p}_{st}, \forall t$ (possibly randomized) with performance guarantee $\tilde{\delta}$, which is explained below: We denote the corresponding virtual queue, Lyapunov drift, and time average revenue expectation by $\tilde{q}_{st}$, $\tilde{\Delta}_{st}$ and $\tilde{v}$, respectively. The $\omega$-only policy satisfies that it is feasible for problem \eqref{eq:v_1}, $\tilde{p}_{st}$ only depends on $\gamma_{st}$ and is independent of the virtual queue for all $t$, and the performance $\tilde{v} \geq v_1^* - \tilde{\delta}$.

By $p_{st}^c p_{st}^d = 0$, $0 \leq p_{st}^c \leq P_s^{max}$, and $0 \leq p_{st}^d \leq P_s^{max}$, there is an upper bound $C_0$ for the quadratic term in $\Delta_{st}$, i.e.,
\begin{align}
    (p_{st}^c \tau \eta_s^c - p_{st}^d \tau / \eta_s^d)^2 / 2 \leq C_0 \triangleq (P_s^{max}\tau/\eta_s^d)^2/2. \nonumber
\end{align}

Because $p_{st}^*$ is optimal in problem \eqref{eq:drift-plus-penalty}, we have
\begin{small}
\begin{align}
    \mathbb{E}[\Delta_{st}^*+V_s(-\gamma_{st} p_{st}^* \tau)|q_{st}^*,\gamma_{st}] \leq \mathbb{E}[\tilde{\Delta}_{st}+V_s(-\gamma_{st} \tilde{p}_{st} \tau)|q_{st}^*,\gamma_{st}], \nonumber
\end{align}
\end{small}
then
\begin{subequations}
\label{eq:omega-only}
\begin{align}
    & \mathbb{E}[\Delta_{st}^*+V_s(-\gamma_{st} p_{st}^* \tau)|q_{st}^*] \\
    \leq~ & \mathbb{E}[\tilde{\Delta}_{st}+V_s(-\gamma_{st} \tilde{p}_{st} \tau)|q_{st}^*] \\
    \label{eq:omega-only-1}
    \leq~ & \mathbb{E}[C_0+ (\tilde{p}_{st}^c \tau \eta_s^c - \tilde{p}_{st}^d \tau /\eta_s^d)q_{st}^* +V_s(-\gamma_{st} \tilde{p}_{st} \tau)|q_{st}^*] \\
    \label{eq:omega-only-2}
    =~ & C_0 + (\tilde{p}_{st}^c \tau \eta_s^c - \tilde{p}_{st}^d \tau /\eta_s^d)q_{st}^* +V_s(-\gamma_{st} \tilde{p}_{st} \tau),
\end{align}
\end{subequations}
where the independence of the $\omega$-only policy on $q_{st}^*$ is utilized from \eqref{eq:omega-only-1} to \eqref{eq:omega-only-2}.

Taking expectations on the two sides of \eqref{eq:omega-only}, we have
\begin{subequations}
\begin{align}
    & \mathbb{E}[\Delta_{st}^*]+V_s\mathbb{E}[-\gamma_{st} p_{st}^* \tau] \\
    \leq~ & C_0 + \mathbb{E}[(\tilde{p}_{st}^c \tau \eta_s^c - \tilde{p}_{st}^d \tau /\eta_s^d)q_{st}^*] +V_s\mathbb{E}[-\gamma_{st} \tilde{p}_{st} \tau] \\
    \label{eq:omega-only-3}
    =~ & C_0 + \mathbb{E}[\tilde{p}_{st}^c \tau \eta_s^c - \tilde{p}_{st}^d \tau /\eta_s^d]\mathbb{E}[q_{st}^*] +V_s\mathbb{E}[-\gamma_{st} \tilde{p}_{st} \tau],
\end{align}
\end{subequations}
where \eqref{eq:omega-only-3} is also because $\tilde{p}_{st}$ only depends on $\gamma_{st}$.

Because $\gamma_{st}, \forall t$ are independent and identically distributed and $\tilde{p}_{st}$ only depends on $\gamma_{st}$, there is a constant $C_1$ so that $\mathbb{E}[\tilde{q}_{s(t+1)}-\tilde{q}_{st}] = \mathbb{E}[\tilde{p}_{st}^c \tau \eta_s^c - \tilde{p}_{st}^d \tau /\eta_s^d] = C_1, \forall t$ under strategy $\tilde{p}_{st}$. By the mean rate stability of the virtual queue,
\begin{align}
    0 & = \lim_{T \rightarrow \infty} \frac{1}{T} \mathbb{E}[\tilde{q}_{sT}] \nonumber \\
    & = \lim_{T \rightarrow \infty} \frac{1}{T} \left(\mathbb{E}[\tilde{q}_{s0}] + \sum_{t = 0}^{T-1} \mathbb{E}[\tilde{q}_{s(t+1)}-\tilde{q}_{st}]\right) \nonumber \\
    & = \lim_{T \rightarrow \infty} \frac{1}{T} (\mathbb{E}[\tilde{q}_{s0}] + T\cdot C_1) = C_1. \nonumber
\end{align}
Thus,
\begin{align}
\label{eq:omega-only-4}
    \mathbb{E}[\Delta_{st}^*]+V_s\mathbb{E}[-\gamma_{st} p_{st}^* \tau] \leq C_0 + V_s\mathbb{E}[(-\gamma_{st} \tilde{p}_{st} \tau)].
\end{align}
Note that $\sum_{t = 0}^{T-1} \mathbb{E}[\Delta_{st}^*] = \mathbb{E}[l_{sT}^* - l_{s0}^*]$ is bounded because $q_{st}^*$ is bounded by Theorem \ref{thm:feasible}. Thus, 
\begin{align}
    \lim_{T\rightarrow \infty} \frac{1}{T} \sum_{t = 0}^{T-1} \mathbb{E}[\Delta_{st}^*] = 0. \nonumber
\end{align}
Then by taking time average in \eqref{eq:omega-only-4}, we have
\begin{footnotesize}
\begin{align}
    & V_s \lim_{T\rightarrow \infty} \frac{1}{T} \sum_{t = 0}^{T-1} \mathbb{E}[-\gamma_{st} p_{st}^* \tau] \leq C_0 + V_s \lim_{T\rightarrow \infty} \frac{1}{T} \sum_{t = 0}^{T-1}\mathbb{E}[(-\gamma_{st} \tilde{p}_{st} \tau)] \nonumber \\
    & \implies - V_s v^* \leq C_0 - V_s \tilde{v} \nonumber \\
    & \implies v^* \geq \tilde{v} - \frac{C_0}{V_s} \geq v_1^* - \tilde{\delta} - \frac{C_0}{V_s}. \nonumber
\end{align}
\end{footnotesize}
Let $\tilde{\delta} \rightarrow 0$, we get \eqref{eq:v_1-inequality}, and the proof is completed.

\setcounter{equation}{0}  
\renewcommand{\theequation}{D.\arabic{equation}}
\renewcommand{\thetheorem}{D.\arabic{theorem}}
\section{Proof of Proposition \ref{prop:bidding}}
\label{appendix-D}

a) We first prove that $f_{st}(p_{st})$ is convex. Because $\lim_{p_{st} \rightarrow 0_-} f_{st}(p_{st}) = \lim_{p_{st} \rightarrow 0_+} f_{st}(p_{st}) = 0$, $f_{st}(p_{st})$ is continuous as a function on $[\underline{P}_{st}, \overline{P}_{st}]$. Its derivative
\begin{align}
    \frac{df_{st}(p_{st})}{dp_{st}} = \left\{
    \begin{array}{ll}
        \frac{p_{st} \tau (\eta_s^c)^2 - q_{st}\eta_s^c }{V_s} - \psi_{st}^*, & \text{if}~ \underline{P}_{st} \leq p_{st} < 0, \\
        \frac{p_{st}\tau - q_{st} \eta_s^d}{V_s (\eta_s^d)^2} - \psi_{st}^*, & \text{if}~ 0 < p_{st} \leq \overline{P}_{st},
    \end{array}
    \right. \nonumber
\end{align}
is nondecreasing and linear on $[\underline{P}_{st},0)$ and $(0,\overline{P}_{st}]$. Moreover,
\begin{align}
    & \frac{df_{st}}{dp_{st}}(0_-) = \lim_{p_{st} \rightarrow 0_-} \frac{f_{st}(p_{st})-f_{st}(0)}{p_{st}} = - \frac{q_{st}\eta_s^c }{V_s} - \psi_{st}^*, \nonumber \\
    & \frac{df_{st}}{dp_{st}}(0_+) = \lim_{p_{st} \rightarrow 0_+} \frac{f_{st}(p_{st})-f_{st}(0)}{p_{st}} = - \frac{q_{st}}{V_s \eta_s^d} - \psi_{st}^*. \nonumber
\end{align}
According to \eqref{eq:V-E-value},
\begin{align}
    \underline{E}_s - E_s = - V_s \overline{\gamma}_s \eta_s^d,~ \overline{E}_s - E_s = -V_s \underline{\gamma}_s / \eta_s^c, \nonumber
\end{align}
so
\begin{align}
    - V_s \overline{\gamma}_s \eta_s^d \leq q_{st} \leq - V_s \underline{\gamma}_s/\eta_s^c \leq 0. \nonumber
\end{align}
Because $q_{st} \leq 0$, $V_s > 0$, and $\eta_s^c, \eta_s^d \in (0,1)$, we have $-q_{st}\eta_s^c /V_s \leq -q_{st}/(V_s \eta_s^d)$.
Thus, $f_{st}(p_{st})$ is convex on $[\underline{P}_{st}, \overline{P}_{st}]$.

b) Then we prove the ES operational constraints \eqref{eq:offline-2}-\eqref{eq:offline-4}, which come down to $p_{st} \in [-P_s^{max}, P_s^{max}]$ and $q_{s(t+1)} \in [\underline{E}_s - E_s, \overline{E}_s - E_s]$. According to \eqref{eq:bounds-ES}:
\begin{itemize}
    \item When $q_{st} \leq - V_s \underline{\gamma}_s / \eta_s^c - P_s^{max} \tau \eta_s^c$, we have $\underline{P}_{st} = - P_s^{max}$ and
    \begin{align}
        q_{s(t+1)} \leq q_{st} + P_s^{max} \tau \eta_s^c \leq - V_s \underline{\gamma}_s / \eta_s^c = \overline{E}_s - E_s. \nonumber
    \end{align}
    \item When $- V_s \underline{\gamma}_s / \eta_s^c - P_s^{max} \tau \eta_s^c \leq q_{st} \leq - V_s \underline{\gamma}_s / \eta_s^c$, 
    \begin{align}
        \underline{P}_{st} = \frac{q_{st} \eta_s^c + V_s \underline{\gamma}_s}{\tau(\eta_s^c)^2} \in [-P_s^{max}, 0], \nonumber
    \end{align}
    and
    \begin{align}
        q_{s(t+1)} \leq q_{st} - \frac{q_{st} \eta_s^c + V_s \underline{\gamma}_s}{\tau(\eta_s^c)^2} \tau \eta_s^c = - \frac{V_s \underline{\gamma}_s}{ \eta_s^c} = \overline{E}_s - E_s. \nonumber
    \end{align}
    \item When $q_{st} \geq - V_s \underline{\gamma}_s / \eta_s^c$, we have $\underline{P}_{st} = 0$ and $q_{s(t+1)} \leq q_{st} \leq \overline{E}_s - E_s$.
\end{itemize}
To sum up, by checking the formula of $\underline{P}_{st}$, we know that $p_{st} \geq -P_s^{max}$ and $q_{s(t+1)} \leq \overline{E}_s - E_s$. Similarly, we have $p_{st} \leq P_s^{max}$ and $q_{s(t+1)} \geq \underline{E}_s - E_s$ from the formula of $\overline{P}_{st}$. Therefore, the charging and discharging power bounds and the SoC bounds all hold under the proposed bidding cost curve $f_{st}(p_{st}), \underline{P}_{st} \leq p_{st} \leq \overline{P}_{st}$.

c) Finally, we prove that the market clearing results coincide with the optimal operation strategy.
The relationship between $p_{st}$ and $\gamma_{st}$ is the operation strategy in \eqref{eq:strategy-gamma}. For clarity, we express it by $p_{st} = h(\gamma_{st})$ using the continuous and nondecreasing function $h$ as follows.
\begin{align}
   & h(\gamma_{st}) \triangleq \nonumber\\
    & \left\{
    \begin{array}{ll}
        - P_s^{max}, & \text{if}~ \gamma_{st} \leq - \frac{q_{st} + P_s^{max} \tau\eta_s^c }{V_s} \eta_s^c, \\
        \frac{q_{st} \eta_s^c + V_s \gamma_{st}}{\tau(\eta_s^c)^2 }, & \text{if}~ - \frac{q_{st} + P_s^{max} \tau\eta_s^c }{V_s} \eta_s^c \leq \gamma_{st} \leq - \frac{q_{st} \eta_s^c}{V_s}, \\
        0, & \text{if}~ - \frac{q_{st} \eta_s^c}{V_s} \leq \gamma_{st} \leq - \frac{q_{st}}{V_s \eta_s^d}, \\
        \frac{q_{st} / \eta_s^d + V_s \gamma_{st}}{\tau / (\eta_s^d)^2}, & \text{if}~ - \frac{q_{st}}{V_s \eta_s^d} \leq \gamma_{st} \leq \frac{P_s^{max} \tau - q_{st} \eta_s^d}{V_s (\eta_s^d)^2}, \\
        P_s^{max}, & \text{if}~ \frac{P_s^{max} \tau - q_{st} \eta_s^d}{V_s (\eta_s^d)^2} \leq \gamma_{st}.
    \end{array}
    \right. \nonumber
\end{align}
We want to prove that $p_{st}^* = h(\gamma_{st}^*)$, where $p_{st}^*$ and $\gamma_{st}^*$ are the power output and combined price in the market clearing results, respectively. 
The notations of other variables are similar. Recall that the market clearing result for ES $s$ is equivalent to solving \eqref{eq:ESindividual} given the LMP $\lambda_{st}^*$. Because $f_{st}(p_{st})$ is convex and differentiable in $(\underline{P}_{st},0)$ and $(0,\overline{P}_{st})$, the optimal solution $p_{st}^*$ satisfies
\begin{small}
\begin{align}
    \left\{
    \begin{array}{ll}
        \frac{df_{st}}{dp_{st}}(p_{st}^*) = \lambda_{st}^*, & \text{if}~ p_{st}^* \in (\underline{P}_{st},0) \cup (0,\overline{P}_{st}), \\
        \frac{df_{st}}{dp_{st}}((\underline{P}_{st})_+) \geq \lambda_{st}^*, & \text{if}~ p_{st}^* = \underline{P}_{st}, \\
        \frac{df_{st}}{dp_{st}}((\overline{P}_{st})_-) \leq \lambda_{st}^*, & \text{if}~ p_{st}^* = \overline{P}_{st}, \\
        \frac{df_{st}}{dp_{st}}(0_-) \leq \lambda_{st}^* \leq \frac{df_{st}}{dp_{st}}(0_+), & \text{if}~ p_{st}^* = 0.
    \end{array}
    \right. \nonumber
\end{align}
\end{small}
We check the four cases one by one.
\begin{itemize}
\item When $p_{st}^* \in (\underline{P}_{st},0)\cup (0,\overline{P}_{st})$,
\begin{align}
    h(\gamma_{st}^*) = h\left(\frac{df_{st}}{dp_{st}}(p_{st}^*) + \psi_{st}^* \right) = h(\gamma(p_{st}^*)) = p_{st}^*. \nonumber
\end{align}
\item When $p_{st}^* = \underline{P}_{st}$, $\lambda_{st}^* \leq \frac{df_{st}}{dp_{st}}((\underline{P}_{st})_+)$. Because $h$ is nondecreasing,
\begin{align}
    h(\gamma_{st}^*) \leq h \left( \frac{df_{st}}{dp_{st}}((\underline{P}_{st})_+) + \psi_{st}^* \right) = \underline{P}_{st}. \nonumber
\end{align}
By $\underline{P}_{st} \leq h(\gamma_{st}^*) \leq \overline{P}_{st}$, we have $h(\gamma_{st}^*) = \underline{P}_{st} = p_{st}^*$.
\item When $p_{st}^* = \overline{P}_{st}$, $\lambda_{st}^* \geq \frac{df_{st}}{dp_{st}}((\overline{P}_{st})_-)$, so
\begin{align}
    h(\gamma_{st}^*) \geq h \left( \frac{df_{st}}{dp_{st}}((\overline{P}_{st})_-) + \psi_{st}^* \right) = \overline{P}_{st}. \nonumber
\end{align}
Then $h(\gamma_{st}^*) = \overline{P}_{st} = p_{st}^*$ follows from $\underline{P}_{st} \leq h(\gamma_{st}^*) \leq \overline{P}_{st}$.
\item When $p_{st}^* = 0$, $\frac{df_{st}}{dp_{st}}(0_-) \leq \lambda_{st}^* \leq \frac{df_{st}}{dp_{st}}(0_+)$.
Then
\begin{align}
    0 & = h\left( \frac{df_{st}}{dp_{st}}(0_-) + \psi_{st}^* \right) \nonumber \\
    & \leq  h(\gamma_{st}^*) \leq h\left( \frac{df_{st}}{dp_{st}}(0_+) + \psi_{st}^*\right) = 0, \nonumber
\end{align}
which implies $h(\gamma_{st}^*) = 0 = p_{st}^*$.
\end{itemize}
This completes the proof.

\bibliographystyle{IEEEtran}
\bibliography{IEEEabrv,mybib}

\end{document}